\documentstyle[12pt,amssymb,amscd]{amsart}

\newtheorem{thm}{Theorem}[section]
\newtheorem{prop}[thm]{Proposition}
\newtheorem{lma}[thm]{Lemma}

\newenvironment{defn}{\trivlist
\item[{\bf\hspace{1.88pt}\addtocounter{thm}{1}
Definition~\arabic{section}.\arabic{thm}.}]}{\endtrivlist}
\newenvironment{remark}{\resultspace\par\noindent{\bf Remark.}}{\resultspace}
 
\def\cstar{$\mathrm{C}^{*}$}

\def\indlimit{{\displaystyle \lim_{\longrightarrow} }}
\def\algindlimit{{\displaystyle \mathop{\mathrm{alg\,lim}}\limits_{\longrightarrow} }}

\def\resultspace{{\vskip 0.075in}}

\def\chisub#1 {{\chi_{\lower2.5pt\hbox{$\scriptstyle #1$}}}}

\def\angb#1{\langle #1 \rangle }

\numberwithin{equation}{section}
\begin{document}

\bibliographystyle{amsalpha209}

\newcommand{\ca}{{\cal A}}
\newcommand{\cb}{{\cal B}}
\newcommand{\cc}{{\cal C}}
\newcommand{\cd}{{\cal D}}
\newcommand{\cf}{{\cal F}}
\newcommand{\cg}{{\cal G}}
\newcommand{\ch}{{\cal H}}
\newcommand{\ck}{{\cal K}}
\newcommand{\cl}{{\cal L}}
\newcommand{\cn}{{\cal N}}
\newcommand{\cm}{{\cal M}}
\newcommand{\cs}{{\cal S}}
\newcommand{\ct}{{\cal T}}
\newcommand{\ot}{\otimes} 
\newcommand{\op}{\oplus}  
\newcommand{\pr}{\prime}


\newcommand{\matrrccc}[6]{\mbox{$ \left[ \begin{array}{ccc}
	    #1&#2&#3\\ #4&#5&#6 \end{array} \right] $}}
\newcommand{\matrccc}[3]{\mbox{$ \left[ \begin{array}{ccc}
	    #1&#2&#3 \end{array} \right] $}}
\newcommand{\matrcccc}[4]{\mbox{$ \left[ \begin{array}{cccc}
	    #1&#2&#3&#4 \end{array} \right] $}}
\newcommand{\matrrrc}[3]{\mbox{$ \left[ \begin{array}{c}
	    #1 \\ #2 \\ #3 \end{array} \right] $}}
\newcommand{\matrrrrc}[4]{\mbox{$ \left[ \begin{array}{c}
	    #1 \\ #2 \\ #3 \\ #4 \end{array} \right] $}}
\newcommand{\matrrrcc}[6]{\mbox{$ \left[ \begin{array}{cc}
	    #1&#2 \\ #3&#4 \\#5&#6 \end{array} \right] $}}
\newcommand{\matrrrrcc}[8]{\mbox{$ \left[ \begin{array}{cc}
	    #1&#2\\#3&#4\\#5&#6\\#7&#8 \end{array} \right] $}}
\newcommand{\matrrrccc}[9]{\mbox{$ \left[ \begin{array}{ccc}                    
            #1&#2&#3\\#4&#5&#6 \\#7&#8&#9 \end{array} \right] $}}

\newcommand{\mattri}[9]{\mbox{$ \left[ \begin{array}{cccc}
   #1\\#2&#3\\#4&#5&#6 \\#7&#8&#9&1 \end{array} \right] $}}
\newcommand{\mtab}[7]{
\begin{array}{rrr rr}\hline\\[-9pt]
  #1	 & &#3 & 0 & 1 \\[1pt]\hline\\[-9pt]
      #2 &0&  &#4 &#5 \\
	 &1&  &#6 &#7 \\[1pt]\hline
\end{array}  }

\def\angb#1{\langle #1 \rangle }

\newcommand{\matrcc}[2]{\mbox{$ \left[ \begin{array}{cc} #1 & #2 \end{array}
           \right] $}}
\newcommand{\matrrc}[2]{\mbox{$ \left( \begin{array}{c}  #1 \\ #2 \end{array}
           \right) $}}
\newcommand{\matrrcc}[4]{\mbox{$ \left[ \begin{array}{cc} #1&#2 \\ #3&#4
           \end{array} \right] $}}
\newcommand{\IZ}{\Bbb{Z}}

\bibliographystyle{amsalpha209}

\title[Relative Positions of Matroid Algebras]{Relative
 Positions of Matroid Algebras}
\author{S. C. Power}
\address{Dept.\ of Mathematics \& Statistics\\ Lancaster University\\
     Lancaster, U.K.\ LA1~4YF}
\email{s.power@@lancaster.ac.uk}
\thanks{Jan 1 1998}
\maketitle
 
\begin{abstract}
A classification is given for regular positions
$D \oplus D \subseteq D$ of Jones index 4 where
$D = \algindlimit M_{n_k}({\Bbb C})$ is an even  matroid algebra and 
where the individual summands have index 2. A similar classification
is obtained for positions of direct sums of
2-symmetric algebras and, in the odd case, for the positions
of sums of 2-symmetric C*-algebras in matroid C*-algebras.
The approach relies on an analysis 
of intermediate non-self-adjoint operator algebras and 
the classifications are given in 
terms of $K_0$ invariants, partial isometry homology
and scales in the composite invariant 
$K_0(-) \oplus H_1(-)$.
\end{abstract}

\section{Introduction}

A matroid algebra is a complex algebra
with involution which is a union of a chain of subalgebras each of
which is a full matrix algebra. Together with their 
C$^{\ast}$-algebra counterparts, the matroid C$^{\ast}$-algebras,
these fundamental algebras are classified in terms of the range of 
a normalised trace
on projections, or, equivalently, by
$K_0$ as a scaled ordered dimension group.
See Dixmier \cite{dix} and Elliott \cite{ell-1}.
Although a direct sum of matroid 
algebras  is similarly determined
we show that the position of
this direct sum as a subalgebra of a matroid superalgebra depends on
a number of much finer invariants.  
Here we say that subalgebras are in the same position if
they are conjugate by an automorphism of the superalgebra.  
The terminology follows the von Neumann algebra usage, as 
in Ocneanu \cite{ocn}.
An analysis is given of the relative position of the summands
of a direct sum subalgebra of a matroid algebra
and in particular
a complete classification is obtained
for what may be viewed as the first nontrivial case, namely
that of regular positions
\[
D \oplus D \subseteq \tilde{D}
\]
where $D$ and $\tilde{D}$ are even matroid algebras
and the individual summand inclusions
\[
D \oplus 0 \subseteq p_1\tilde{D}p_1, \ \ \  0 \oplus D 
\subseteq p_2\tilde{D}p_2
\]
for the central projections $p_1, p_2$, 
are standard inclusions of Jones index 2.
Thus $D \oplus D$ has index 4 in $\tilde{D}$. See 
\cite{goo-del-jon} and 
\cite{wat}. The regularity of the inclusion
requires that the subalgebra and the superalgebra share a 
standard regular AF masa.

Classifications are given in terms of
\medskip

(a) the scaled ordered group 
$(K_0(-), \Sigma_0(-))$ of the subalgebra,
\medskip
 
(b) a scaled partial isometry homology group $(H_1(-), \Sigma_1(-))$, 
which may be realised 
as a subgroup of ${\Bbb Q}$ together 
with an interval symmetric about the origin,
\medskip   
 
(c) a joint scale $\Sigma(-)$ in the composite invariant 
$K_0(-) \oplus H_1(-)$.
\medskip

\noindent

As the index 2 positions are unique the invariants 
which augment $K_0$ 
should be viewed as invariants for the
relative position of the matroid algebra summands.

The joint scale in the composite invariant, which need 
not coincide with the product scale,
accounts for a number of obstructions to lifting 
$K_0 \oplus H_1$ isomorphisms
to algebra isomorphisms.
The simplest of these  is an $H_0H_1$ coupling obstruction
which manifests itself as a coset in 
${\Bbb R/Q}$.
If this vanishes then so too
do the other obstructions
and the partial isometry homology group $H_1(-)$
is the remaining invariant for relative position.
Thus  there are uncountably many
positions in this case and there is a unique
position with trivial $H_1$ invariant.

Similar classifications are obtained for
inclusions of direct sums of the so called 2-symmetric
algebras in matroid algebras. In the odd case, which is somewhat
more accessible,
we  obtain an approximate version of the key arguments of section 3
and this enables
the  classification of inclusions of sums of
2-symmetric C*-algebras.
 
If the partial isometry homology invariant
 $H_1(\pi)$ is nonzero for the ordered inclusion
$\pi : D \oplus D \to D$ then we calculate the outer
automorphism group $Out_{D \oplus D}(D)$
of 
\\ 
subalgebra-respecting automorphisms modulo
automorphisms with approximately inner restrictions.
In the unital matroid case this group is equal to ${\Bbb Z}_2$
if the $H_0H_1$ coupling obstruction is nontrivial and 
coincides with $Aut(H_1(\pi))$
otherwise. The generator for  ${\Bbb Z}_2$ derives from  a
homology inverting automorphism.
On the other hand, in the 2-symmetric case we show that there 
can be intriguing 
obstructions to homology inversion.
This is the case, for example, for positions determined
by the so called homologically asymmetric 4-cycle algebra systems 
identified in Donsig and Power \cite{don-pow-2}.
The determination of $Out_{D \oplus D}(D)$
when  $H_1(\pi)$ is trivial is left open.

If a partly self-adjoint algebra $A$ is a generating subalgebra 
of a self-adjoint algebra $B$ then the position
$A \cap A^* \subseteq B$ is an invariant for the 
star extendible isomorphism class of $A$.
Dually, the isomorphism type of the algebras in the intermediate
subalgebra lattice of an inclusion provide invariants for that inclusion.
This natural link between self-adjoint 
and non-self-adjoint algebra is
essential to our approach
and it is by considering direct limits  of
non-self-adjoint finite-dimensional algebras that 
we obtain classifying invariants and diverse inclusions.

For a regular direct system  of general digraph
algebras (finite-dimensional incidence algebras) such as
\[
A_1 @>{\alpha_1}>> A_2 @>{\alpha_2}>>
        A_3 \to \dots,
\]
the natural
limit homology groups 
\[
H_n(\{A_k, \alpha_k \}) = 
\indlimit(H_n(A_k), H_n (\alpha_k))
\]
were introduced in Davidson and Power
\cite{dav-pow} and given intrinsic formulations, as partial isometry
homology, in \cite{pow-hom2},  \cite{pow-hom3}.  
It is necessary here to consider
regular morphisms, that is,  homomorphisms 
that  are direct sums of multiplicity
one embeddings  (\cite{pow-ko}), and in this case 
the limit homology groups are well-defined invariants for 
the regular isomorphism of regular direct systems.
With this invariant and the introduction of new 
scales  for $K_0 \oplus H_1$
we completed the 
$K_0H_1$ classification
of certain direct systems of 4-cycle algebras
up to  regular isomorphism. 
See Donsig and Power \cite{don-pow-2}.
On the other hand in \cite{don-pow-1} it was shown 
that such direct systems 
can be irregularly isomorphic, in a nontrivial manner,
and so the classification
of even the algebraic direct limits of the systems was left open.
In the present paper we overcome this obstacle by establishing
the well-definedness of the 
partial isometry homology group
invariants. This 
also leads to the correct homology invariants for
the discrimination and classification of subalgebra positions.

The paper is organised as follows. 
In Section 2 we define the standard index 2 inclusions of matroid algebras,
the generic index 4 inclusions of direct sums of matroid 
algebras, inclusions of 2-symmetric algebras, and
generic limits of 4-cycle algebras.  
The  intrinsic formulation of $H_n(\{A_k, \alpha_k\})$ 
as the homology of a chain complex is also
indicated.
In Section 3 we obtain the key fact
that irregular isomorphisms  of regular 
systems  of 4-cycle algebras
are only possible if the systems are of 
matroid type and have vanishing homology.
In sections 4 and 5 we classify generic limits of 4-cycle algebras,
relative positions of direct sums of matroid algebras and 2-symmetric algebras,
and we determine
$Out_{D_1 \oplus D_2}(D)$. We also indicate obstructions to homology inverting
automorphisms.
 These sections 
are largely algebraic and concern only algebraic direct limits.

In section 6 
the key arguments of sections 3 and 4 are generalised
to the case of approximately regular embeddings and the classification
of operator algebra limits of generic 4-cycle algebra systems is
obtained, in 
the odd case, in terms of
\[
(K_0(-) \oplus H_1(-), \Sigma_0(-), \Sigma(-)).
\]
The stability of $K_0H_1$ invariants obtained here
leads also to the fact that close limit algebras are isomorphic,
in analogy with situation for AF C*-algebras \cite{chr}, \cite{phi-rae}.

The reader may notice a distinct parallel between the $K_0H_1$
classification of cycle algebras and the modern
$K_0K_1$ determination of C*-algebras, such as is indicated in
Elliott \cite{ell-2} and Riordam \cite{rio}. 
Thus homology inverting automorphisms
bear an analogy with $K_1$ inverting automorphisms,
such as the flip on the Bunce-Deddens algebras.
The results below suggest extending the
well-definedness of partial isometry homology further in
order to approach the $K_*H_*$ determination of partly
self-adjoint subalgebras of such C*-algebras.

I would like to thank Alan Donsig for some helpful discussions.

\section{Generic Inclusions}

Let $F_{ {q}}$ be an (unclosed) unital matroid 
algebra associated with the generalised
integer $ {q} =q_1q_2 \dots$, with $q_i \geq 2$ for all $i$.  Suppose that
$ {q}=2^{\infty}  {p}$ with $ {p} = p_1p_2 \dots$, and  $p_i \geq 2$ for
all $i$ and consider the realisation of $F_{q}$ as a direct limit

\[
F_{{q}} = \algindlimit (M_{2^k{p_1 \dots p_k}} \oplus M_{2^k{p_1
\dots p_k}}, \alpha_k)
\]

\noindent where 

$$ \begin{array}{l}
K_0 \alpha_k =
\matrrcc{p_{k+1}}{p_{k+1}}{p_{k+1}}{p_{k+1}} \end{array}$$

\noindent for all $k \geq 1$.  
View the maps $\alpha_k$ as the restrictions to the
block diagonal subalgebras of the map 

$$\alpha_k : M_{2^{k+1}p_1 \dots p_k} \rightarrow 
M_{2^{k+2}p_1 \dots p_{k+1}}$$

\noindent for which
\medskip

\[
\alpha_k : \matrrcc{a}{b}{d}{c} \quad \rightarrow \quad \quad   
\left[ \begin{array}{cccc}
\sigma_k(a) &0&0&\sigma_k(b)\\
0&\sigma_k(c)&\sigma_k(d)&0\\
0&\sigma_k(b)&\sigma_k(a)&0\\ 
\sigma_k(d) &0&0&\sigma_k(c)
\end{array} \right]
\]
\medskip

\noindent where $\sigma_k (a) = a \oplus \dots \oplus 
a$ \ \ ($p_{k+1}$ times).  

In this way  obtain an inclusion $F_{ {q}} \rightarrow F_{2{q}} \ \ (=
F_{ {q}})$ with \ $F_{2{q}} = F_{{q}}+ JF_{{q}}$ where $J$ is the self-adjoint
unitary determined by the symmetry 
$\begin{array}{l}
\matrrcc{O}{I}{I}{0} \end{array}$
in $M_{4{p_1}}$.  For convenience we
refer to this inclusion $F_{{q}} \rightarrow F_{2{q}}$ as the
{\it standard regular inclusion } of index 2 and we remark that it 
also arises from the crossed product $F_{2q} = F_q \times_\alpha {\Bbb Z}_2$
where $\alpha$ is a product type symmetry. (See \cite{fac-mar}.
Note that  the abelian diagonal subalgebra

$$C = \algindlimit (C_k, \alpha_k)$$

\noindent in $F_{2{q}}$ is a {\it regular} diagonal
masa for both $F_{2{q}}$ and
$F_{{q}}$. That is, each of the algebras $F_{2{q}}$ and
$F_{{q}}$ is generated by the
partial isometries that normalise $C$.

If $F_{ {q}} \oplus F_{ {q}} \rightarrow F_{2  {q}} \oplus F_{2
 {q}}$ is the associated index 2 position then it may seem
curious, given the uniqueness of the trivial
index 2 position
$F_{2  {q}} \oplus F_{2  {q}} \rightarrow F_{2  {q}} \otimes M_2 $,
that there are uncountably many positions. Nevertheless this 
is a consequence of Theorem 4.5.
\medskip

We now consider embeddings of
non-self-adjoint 
algebras 
which lead to
an apparent variety of
relative positions for matroid algebras
and 2-symmetric algebras.
\medskip 

Let $A \subseteq M_n$ be a complex algebra which contains a maximal
abelian self-adjoint subalgebra (masa).  Then there is a  matrix unit
system 
$\{e_{ij} : 1 \leq i, \  j \leq n \}$ for $M_n$ such that
the masa is spanned by the diagonal matrix units
and $A$ itself is spanned
by matrix units.
We refer to  $A$ as a {\it digraph algebra} whose digraph is the
directed graph on $n$ vertices with edges $(i,j)$ directed from $j$ to
$i$, for each $e_{ij}$ in $A$.  Let $\phi : A_1 \rightarrow A_2$ be
an embedding of digraph algebras with $\phi(C_1) \subseteq C_2$ for some
masas $C_i \subseteq A_i$, $i=1,2$, and assume that $\phi$ is
{\it star extendible} in the sense that there is a (necessarily unique)
C$^{\ast}$-algebra extension ${\phi}: C^{\ast}(A_1) \rightarrow
C^{\ast} (A_2)$.  Then the following statements are equivalent.

\medskip

(a) $\phi$ is a direct sum of multiplicity one embeddings.
\medskip

(b) $\phi$ maps the (partial isometry) normaliser of $C_1$ into the
normaliser of $C_2$.
\medskip

In this event we say that $\phi$ is a {\it regular}
star extendible embedding.  
\medskip
  
A 4-cycle matrix algebra is a complex subalgebra of $M_r$
consisting of
partitioned matrices of the form
\[
 \left[ \begin{array}{cccc}
a_{11}&0&a_{13}&a_{14}\\
0&a_{22}&a_{23}&a_{24}\\
0&0&a_{33}&{0}\\
0&0&0&a_{44}
\end{array} \right]
\]
where $a_{ij}$ is a matrix in $M_{n_i,n_j}$,
the space of $ n_i$ by $n_j$ matrices, and
$r = n_1 + n_2 + n_3 + n_4$.
The simplest of these algebras are
the 4-cycle algebra $A_1$ in $M_4$, which is also denoted $A(D_4)$, and
the equidimensional 4-cycle algebras, which are the algebras $A_1 \otimes M_n$.
 Along with their 2$n$-cycle counterparts,
these algebras have featured prominently in the consideration of 
Hochschild cohomology for non-self-adjoint
operator algebras.  In fact 
\[
\mbox{Hoch}^1(A_1 \otimes
M_n; A_1 \otimes M_n) = {\Bbb C}
\]
for all $n$. See \cite{gil-smi-1},
\cite{dav-pow}, for example.
 
Taking a more geometric
homological perspective, the algebra $A_1 \otimes M_n$ 
determines, through
a matrix unit system, a simplicial complex $\triangle (A_1 \otimes
M_n)$ and therefore one may associate the integral 
simplicial homology group, $H_1 (\triangle
(A_1 \otimes M_n)) = {\Bbb Z}$.  This simplicial complex is the direct product
of the complete complex on $n$ vertices with the 4 cycle complex
$\triangle (A_1)$ indicated by the diagram 
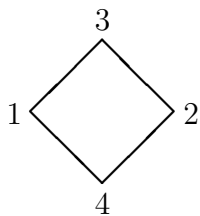
\begin{center}
\setlength{\unitlength}{0.0125in}%
\begin{picture}(110,96)(65,718)
\thicklines
\put(120,790){\line(-1,-1){30}}
\put( 90,760){\line( 1,-1){30}}
\put(120,730){\line( 1, 1){30}}
\put(150,760){\line(-1, 1){30}}
\put( 80,755){\makebox(0,0)[lb]{\raisebox{0pt}[0pt][0pt]{1}}}
\put(154,755){\makebox(0,0)[lb]{\raisebox{0pt}[0pt][0pt]{2}}}
\put(117,794){\makebox(0,0)[lb]{\raisebox{0pt}[0pt][0pt]{3}}}
\put(117,717){\makebox(0,0)[lb]{\raisebox{0pt}[0pt][0pt]{4}}}
\end{picture}
\end{center}
 
More generally, this homology group association is available for
any  digraph
algebra; edges of the digraph $G$ of $A$ contribute 1-simplices
to $\Delta(A)$ and complete (undirected)
subgraphs of $G$ on $t+1$ vertices determine the $t-$simplices.
Note that regular homomorphisms induce homology group homomorphisms.
 
An alternative, homologically equivalent association is to take $\Delta(A)$
to be the complex
determined by the reduced digraph of the digraph algebra $A$.
(Then $\Delta(A_1 \otimes M_n) = \Delta(A_1)$.)

Consider now the 4-cycle digraph algebra $A_1 = A(D_4)$.
Its digraph $D_4$ has edges (1,3), (1,4),
(2,4), (2,3), (1,1), (2,2), (3,3), (4,4).  A multiplicity one
star-extendible embedding
\[
\phi: A(D_4) \otimes M_n \rightarrow A(D_4) \otimes M_m
\]
is said to be {\it rigid} if the images of the rank one partial isometries
$ e_{13} \otimes p, \; e_{14} \otimes p, \; \\ 
e_{24} \otimes p, \;
e_{23} \otimes p,
$
with $p$ a rank one projection of $M_n$, 
are inequivalent rank one
partial isometries in the block off-diagonal part of the range algebra. 
Since the initial and final projections
of these images are rank one projections in
$A(D_4) \otimes M_m$ it follows that 
the images belong, in some order,  to
the distinct block subspaces 
\[
e_{13} \otimes M_m, 
e_{14} \otimes M_m, \\  e_{24} \otimes M_m, e_{23} \otimes M_m.
\]
There are precisely four inner equivalence classes of such embeddings
and these classes correspond to the four symmetries of $D_4$.  
\medskip

\begin{defn}
An
 embedding $\alpha$ between 4-cycle algebras is said to be
{\it generic} if it
is conjugate (inner unitarily equivalent)
to a direct sum of rigid rank one embeddings, with at
least one summand of each of the four types.  
\end{defn}
\medskip

A general direct sum of multiplicity one rigid embeddings,
without the 
constraint of the definition,  is also 
called rigid. Generic and rigid embeddings
between general (nonequidimensional) 4-cycle algebras are defined
similarly.
\medskip

\begin{defn} 
(i) A 4-cycle algebra direct system $\ca = (A_k, \alpha_k)$
is said to be a generic
if infinitely many of the embeddings $\alpha_k$ are generic.
In this case  the associated limit algebra $A =
\displaystyle \algindlimit (A_k, \alpha_k)$
is said to be a generic.

 \ \ \ \ \  \ \ \  (ii) 
An inclusion $D_1 \oplus D_2 \rightarrow B$
of matroid algebras 
is said to be a generic index 4 regular inclusion if it is conjugate,
by a star automorphism of $B$, to an inclusion
$A \cap (A)^* \to B$
where $A$ is a generic 4-cycle limit algebra
as in (i), and $B$ is the 
self-adjoint algebra generated by $A$. In this case the 
system $\{A_k, \alpha_k\}$
is said to be of matroid type.
\end{defn}
\medskip
 
If $\ca$ is a generic system then it has a subsystem in which all 
the embeddings are generic. Henceforth  we shall assume this 
normalisation throughout the paper without further comment.

For an explicit distinguished example consider the 
stationary system $\{A_k, \alpha_k\}$
where  $A_k =  A(D_4) \otimes M_{4^k}$ and
$\alpha_k$ is the embedding
\medskip

\[
 \left[ \begin{array}{cccc}
a&&x&y\\
&b&w&z\\
& &c& \\
 & & &d \end{array} \right]
\to
\left[ \begin{array}{cccccccccccccccc}
a&0&0&0 &&&&  &x&0&0&0&0&0&0&y\\
0&a&0&0 &&&&  &0&0&0&y&x&0&0&0\\
0&0&b&0 &&&&  &0&0&z&0&0&w&0&0\\
0&0&0&b &&&&  &0&w&0&0&0&0&z&0\\
&&&& a&0&0&0  &0&0&y&0&0&x&0&0\\
&&&& 0&a&0&0  &0&x&0&0&0&0&y&0\\
&&&& 0&0&b&0  &w&0&0&0&0&0&0&z\\
&&&& 0&0&0&b  &0&0&0&z&w&0&0&0\\
&&&& &&&&  c&0&0&0 &&&&\\
&&&& &&&&  0&c&0&0 &&&&\\
&&&& &&&&  0&0&d&0 &&&&\\
&&&& &&&&  0&0&0&d &&&&\\
&&&& &&&&  &&&& c&0&0&0\\
&&&& &&&&  &&&& 0&c&0&0\\
&&&& &&&&  &&&& 0&0&d&0\\
&&&& &&&&  &&&& 0&0&0&d
\end{array} \right]
\]
\medskip

One can check that a rigid system is not generic if and only if 
for all large $k$ each of the embeddings $\alpha_k$ is a sum
of multiplity one rigid embeddings
corresponding to at most
two classes.

Not every inclusion $A \cap A^* \to B$ arising from a generic
limit algebra $A$ as in (i)
is a direct sum of matroid algebras. 
In general $A \cap A^*$ is a  direct sum of
two 2-symmetric algebras and we also refer to these inclusions as
generic inclusions (of unspecified index).
\medskip

\begin{defn}
A {\it 2-symmetric algebra} is  the algebraic
direct limit algebra of a system of algebras
each of which is a direct sum of two matrix algebras
and where the 
partial embeddings have 
multiplicities indicated by a symmetric
2 $\times $2 integral matrix of the form
$\bigl[ \begin{smallmatrix} x & y \\ y & x \end{smallmatrix} \bigr]$.
\end{defn}
\medskip

The explicit $K_0$ classification of these algebras
and their C*-algebra closures, the 2-symmetric C*-algebras,
is given in Fack and Marechal \cite{fac-mar} for
the unital equidimensional case,
and in Donsig and Power \cite{don-pow-2} for the general case.
\medskip

The most conspicuous invariant for the subalgebra position
$D_1 \oplus D_2 \subseteq B$ is the lattice of intermediate algebras 
partially ordered by inclusion.
That there are thirteen such algebras follows readily from the fact
that the intermediate algebras are generated by matrix units.
(See also the discussion of inductivity in Chapter 4 of \cite{pow-book}.)
The intermediate subalgebra lattice is indicated in the diagram below 
in the case of a generic unital matroid algebra inclusion
with $D_1 = D_2$
where 
$D \subseteq \tilde{D}$ indicates the standard index 2 position.
The properly nonselfadjoint intermediate algebras
are the algebras $A, A^*$ together with
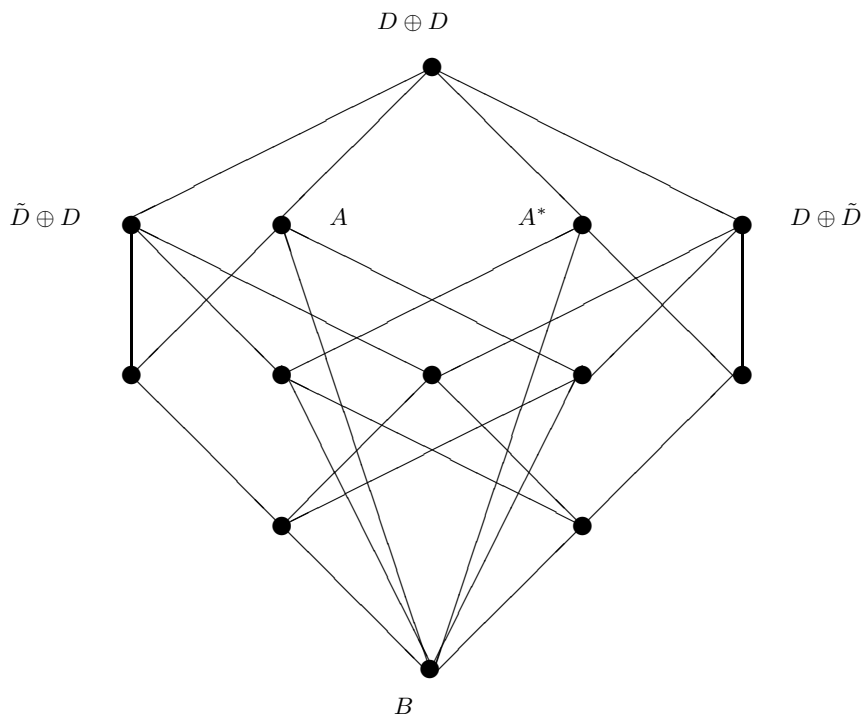
\[
 \langle A,\tilde{D} \oplus D \rangle, \langle A^*,\tilde{D} \oplus D\rangle,
 \langle A, D \oplus \tilde{D}\rangle, 
 \langle A^*, D \oplus \tilde{D}\rangle,
\langle A,\tilde{D} \oplus \tilde{D}\rangle, \langle A^*,
\tilde{D} \oplus \tilde{D}\rangle.
\]


\begin{center}
\setlength{\unitlength}{0.01250000in}%
\begingroup\makeatletter\ifx\SetFigFont\undefined
\def\x#1#2#3#4#5#6#7\relax{\def\x{#1#2#3#4#5#6}}%
\expandafter\x\fmtname xxxxxx\relax \def\y{splain}%
\ifx\x\y   
\gdef\SetFigFont#1#2#3{%
  \ifnum #1<17\tiny\else \ifnum #1<20\small\else
  \ifnum #1<24\normalsize\else \ifnum #1<29\large\else
  \ifnum #1<34\Large\else \ifnum #1<41\LARGE\else
     \huge\fi\fi\fi\fi\fi\fi
  \csname #3\endcsname}%
\else
\gdef\SetFigFont#1#2#3{\begingroup
  \count@#1\relax \ifnum 25<\count@\count@25\fi
  \def\x{\endgroup\@setsize\SetFigFont{#2pt}}%
  \expandafter\x
    \csname \romannumeral\the\count@ pt\expandafter\endcsname
    \csname @\romannumeral\the\count@ pt\endcsname
  \csname #3\endcsname}%
\fi
\fi\endgroup
\begin{picture}(369,299)(235,413)
\thinlines
\put(475,620){\circle*{8}}
\put(349,620){\circle*{8}}
\put(286,620){\circle*{8}}
\put(542,620){\circle*{8}}
\put(475,557){\circle*{8}}
\put(349,557){\circle*{8}}
\put(286,557){\circle*{8}}
\put(542,557){\circle*{8}}
\put(412,686){\circle*{8}}
\put(412,557){\circle*{8}}
\put(349,494){\circle*{8}}
\put(475,494){\circle*{8}}
\put(411,434){\circle*{8}}
\put(412,686){\line(-1,-1){ 64.500}}
\put(349,620){\line(-1,-1){ 63}}
\put(286,557){\line( 0, 1){ 63}}
\put(286,620){\line( 1,-1){ 63}}
\put(349,557){\line( 2, 1){126}}
\put(475,620){\line( 1,-1){ 65}}
\put(542,557){\line( 0, 1){ 63}}
\put(542,620){\line(-1,-1){ 65}}
\put(475,557){\line(-2, 1){126}}
\put(412,686){\line( 1,-1){ 64.500}}
\put(542,620){\line(-2,-1){129.200}}
\put(412,557){\line(-2, 1){126}}
\put(286,557){\line( 1,-1){ 63}}
\put(349,494){\line( 1,-1){ 63}}
\put(412,431){\line( 1, 1){ 63}}
\put(475,494){\line( 1, 1){ 65}}
\put(475,557){\line(-2,-1){126}}
\put(349,494){\line( 1, 1){ 63}}
\put(412,557){\line( 1,-1){ 63}}
\put(475,494){\line(-2, 1){126}}
\put(412,686){\line(-2,-1){127.200}}
\put(412,686){\line( 2,-1){130.400}}
\put(475,620){\line(-1,-3){ 63}}
\put(412,431){\line(-1, 3){ 63}}
\put(349,561){\line( 1,-2){ 63}}
\put(412,435){\line( 1, 2){ 63}}
\put(369,620){\makebox(0,0)[lb]{\smash{\SetFigFont{9}{10.8}{rm}
$A$
}}}
\put(448,620){\makebox(0,0)[lb]{\smash{\SetFigFont{9}{10.8}{rm}
$A^*$
}}}
\put(389,702){\makebox(0,0)[lb]{\smash{\SetFigFont{9}{10.8}{rm}
$D \oplus   D$
}}}
\put(562,620){\makebox(0,0)[lb]{\smash{\SetFigFont{9}{10.8}{rm}
$D \oplus   \tilde{D}$
}}}
\put(235,620){\makebox(0,0)[lb]{\smash{\SetFigFont{9}{10.8}{rm}
$\tilde{D}\oplus D$
}}}
\put(396,415){\makebox(0,0)[lb]{\smash{\SetFigFont{9}{10.8}{rm}
$B$
}}}
\end{picture}
\end{center}
 \medskip

To identify finer invariants for the inclusion
we may focus on the pair $\{A, A^*\}$.
 \medskip

For a rigid embedding $\alpha$ between 4-cycle algebras define the
{\it multiplicity signature} of $\alpha$
as the ordered 4-tuple 
$\{r_1, r_2, r_3, r_4\}$ where $r_1, \dots ,r_4$ are the 
number of summands of each of the four types.
If $r_2$ and $r_4$ 
correspond to the number of rank one rigid embeddings associated with the two
reflection symmetries of $D_4$ then 
the induced $H_1$ group homomorphism $H_1(\alpha)$ 
can be identified with the map
$ {\Bbb Z} \to  {\Bbb Z}$
which is  multiplication by
$r_1 - r_2 + r_3 - r_4.$
This group homomorphism is determined up to sign,
this sign being fixed
 only after
an identification of the reduced graphs of $A_1$ and $A_2$ and realisations 
of $H_1(A_1)$ and $H_1(A_2)$.

It is an elementary but significant fact that the inner equivalence class
of $\alpha$ is determined by the pair 
$K_0(\alpha), H_1(\alpha)$.
See \cite{pow-book}, \cite{don-pow-2}.

More generally one may define the induced homology map
$H_n(\alpha)$ for any regular embedding between digraph algebras, and for 
a regular star extendible direct system
$\{A_k, \alpha_k\}$ one has the abelian group

\[
H_n(\{A_k, \alpha_k\}) = \indlimit (H_n(A_k), H(\alpha_k)).
\]
This is evidently  a well-defined invariant for the 
regular isomorphism of systems. For the stationary example above 
it is the zero group for all $n \ge 1$.

One may define $H_n(\{A_k, \alpha_k\})$ more intrinsically 
in the following manner. For related forms and variations see
Power \cite{pow-crelle}, \cite{pow-hom2}, \cite{pow-hom3}.

Let $Pisom_{reg}(\{A_k, \alpha_k\})$ be the set of partial isometries
$v$ in $A_k$, for some $k$, which are {\it regular}.
By this we mean that $v$ is a sum of rank one partial isometries $w$
for which both $w^*w$ and $ww^*$ belong to $A_k$. The unitary equivalence 
classes of such partial isometries in a fixed algebra $A_k$ are 
in bijective correspondence 
with the edges of the digraph of $A_k$. Define the 1-chain group of 
$\{A_k, \alpha_k\}$ to be the free group generated by all these unitary
equivalence classes $[v]$, for all $k$, subject to the relations
arising from the regular embeddings $\alpha_k$. 
Thus $[v] = [v_1] + [v_2] + \dots + [v_m]$ if $v$ in $A_k$ is equal
to $v_1 + v_2 + \dots + v_m$ in $A_l$ for $l > k$.
Whereas partial
isometries generate 1-simplices, triangles 
of partial
isometries in $\{A_k, \alpha_k\}$ 
(with $v_3 = v_1v_2$  and $v_1^*v_1 = v_2v_2^*$)
determine 2-simplices and a
corresponding 2-chain group  (again, modulo inclusion relations).
There are obvious boundary maps and $H_1(\{A_k, \alpha_k\})$ is defined to be
the appropriate homology group.

The group $H_0(\{A_k, \alpha_k\})$ can also be defined as
the 0-chain group of $\{A_k, \alpha_k\}$ which is similarly defined in terms of
projection classes $[p]$ and it is routine to verify that
\[
H_0(\{A_k, \alpha_k\}) = \indlimit (K_0(C^*(A_k)), K_0(\tilde{\alpha_k})),
\]
where $\tilde{\alpha}_k$ is the star extension of ${\alpha}_k$.
\medskip

We have introduced the generic index 4 inclusions
of matroid algebras by means of systems
of 4-cycle algebras with generic embeddings and this 
presentation emphasises the strong connection with
non-self-adjoint algebras which will be important in the sequel..
The following proposition provides an alternative
direct definition.

First we note that there is a distinguished
index 4 inclusion arising from a rigid nongeneric system.
In the unital equidimensional case this is the natural   
inclusion of $F_q \otimes {\Bbb C}^2$ in $(F_q + JF_q) \otimes M_2$
and we refer to this as the homologically extreme inclusion.

\begin{prop} Let $D_1 \oplus D_2$ be a unital subalgebra of $D_3$
where $D_1, D_2$ and $D_3$ are isomorphic to the unital even matroid
algebra $F_q$.
Then the position  $D_1 \oplus D_2 \subseteq D_3$ is a generic index 4
inclusion (or the homologically extreme inclusion)
if and only if there is a regular presentation
\[
D_3 = \algindlimit (M_{2^ks_k} \otimes M_4, \alpha_k)
\]
\[
D_1 \oplus D_2 = \algindlimit (M_{2^ks_k} \otimes {\Bbb C}^4, \beta_k)
\]
for which
 
(i)
$\beta_k$ is the restriction of $\alpha_k$ to the block diagonal
and has the form
\[
\beta_k(a \oplus b \oplus c \oplus d) = (\sigma_k(a) \oplus
\sigma_k(b)) \oplus (\sigma_k(a) \oplus
\sigma_k(b)) \oplus (\sigma_k(c) \oplus
\sigma_k(d)) \oplus (\sigma_k(c) \oplus
\sigma_k(d)) 
\]
where $\sigma_k$ is the standard embedding of multiplicity
$q_{k+1}$,

(ii) if $p_1$ and $p_2 $ are the orthogonal central projections
of $D_1 \oplus D_2$ then the inclusions
$D_i \subseteq p_iD_3p_i, i = 1, 2$ are standard index 2 inclusions.
\end{prop}

\section{Regular and Irregular Factorisations}

The next lemma 
is crucial and will be 
used with Lemma 3.2 to obtain the well-definedness
of partial isometry homology for algebraic direct limits.
\medskip

The following definitions will be convenient.
Let $\phi : A_1 \rightarrow
A_2 $  be a star extendible embedding of 4-cycle algebras.
Then $\phi$ is said to be {\it locally regular} if the image of each rank 
one partial isometry in $A_1$ is a regular partial isometry
in $A_2$, that is, each such image is an orthogonal sum of rank one
partial isometries in $A_2$ whose
initial and final projections lie in $A_2$.
Also $\phi$ is said to be a {\it proper} embedding if, firstly,
$\phi(A_1^r) \subseteq A_2^r$, where $A_i^r$ denotes the block upper 
triangular part of $A_i$, and secondly, for each rank one partial
isometry $v$ in $A_1$  the partial isometry $\phi(v)$ has
support in each of the four block subspaces of $A_2^r$ . Equivalently
put, 
all the matrices for $\phi$ given in Table 1 are nonzero. The generic 
embeddings defined earlier 
are therefore the same as the proper rigid embeddings.

\begin{lma} { \bf Factorisation dichotomy.}\ \    Let $\phi : A_1 \rightarrow
A_2,$ $\psi: A_2 \rightarrow A_3$,
 $\eta : A_3 \rightarrow
A_4,$ 
be star extendible proper
embeddings between
4-cycle algebras. If the compositions
$\psi \circ \phi$,
$\eta \circ \psi$,
are rigid embeddings then either $\eta$ is locally regular or 
$H_1(\psi \circ \phi) = 0$.
\end{lma}
\medskip

Before beginning the proof we establish some notation for 
the triple $\phi, \psi, \psi \circ \phi$.
\medskip

Assume that $A_1 = A (D_4)$, that $A_2$ and $A_3$ are equidimensional
and that the maps are all unital.   This is the essential case to consider.
\medskip

Choose matrix units for $A_2 \cap A_2^{\ast}$ so that $\phi | (A_1 \cap 
A_1^{\ast})$ has the form
\medskip

$$
\phi : \lambda_1 \oplus \lambda_2 \oplus \lambda_3 \oplus \lambda_4 
\rightarrow (\lambda_1 P' + \lambda_2 {R}') \oplus (\lambda_1 Q' +
\lambda_2 S') \oplus (\lambda_3 P + \lambda_4 R) \oplus (\lambda_3 Q +
\lambda_4 S)
$$
\medskip

\noindent
where $P', R', Q', S', P, R, Q, S$ are orthogonal
projections which are sums of diagonal matrix units.  Set $p'=$
rank $P{'}$, $r'=$ rank $R', \dots, s=$ rank $S$.  Since $\phi$ is
star extendible we have
\medskip
 
\noindent
\[
p'+ q' = r'+s' = p+q = r+s = \rho,
\]
\medskip
 
\noindent
 where $\rho$
is the multiplicity of $\phi$.  Furthermore
\medskip
 
\noindent
\[
\phi : e_{13} \to 
  \begin{bmatrix} 0 & 0 & & & \alpha_1 & 0 & \beta_1 & 0 \\ 0 & 0 & & & 0 & 0 & 0 & 0 \\
	& & 0 & 0 & \delta_1 & 0 & \gamma_1 & 0 \\ & & 0 & 0 & 0 & 0 & 0 & 0 \\ 
	& & & & 0 & 0 & & \\ & & & & 0 & 0 & & \\ & & & & & & 0 & 0 \\
	& & & & & & 0 & 0 \end{bmatrix}
\]
\medskip
 
\noindent
where  { $\matrrcc{\alpha_1}{\beta_1}{\delta_1}{\gamma_1} $}
is a unitary
matrix in  $M_\rho$.  Similarly $\phi (e_{14}), \phi(e_{24}), \phi
(e_{23}) $ have associated unitary matrices
\medskip
 
\noindent
\[
\matrrcc{\alpha_2}{\beta_2}{\delta_2}{\gamma_2},\ \ 
\matrrcc{\alpha_3}{\beta_3}{\delta_3}{\gamma_3},\ \ 
\matrrcc{\alpha_4}{\beta_4}{\delta_4}{\gamma_4}.\ \ 
\]
\medskip
 
\noindent
The dimensions of the matrix entries of these matrices are indicated
in the following table
\medskip
 
\noindent

\begin{center}
\begin{tabular}{|l|l|l|l|l|}
\hline  &&&&\\
& $p$ & $r$ & $q$ & $s$ \\
\hline  &&&&\\
$p'$ & $\alpha_1$ & $\alpha_2$ & $\beta_1$ & $\beta_2$ \\
\hline &&&&\\
$r'$ ~~~~~ & $\alpha_4$ ~~~~~ & $\alpha_3$ ~~~~~ & $\beta_4$ ~~~~~ & $\beta_3$ \\
\hline &&&&\\
$q'$ & $\delta_1$ & $\delta_2$ & $\gamma_1$ & $\gamma_2$ \\
\hline  &&&&\\
$s'$ & $\delta_4$ & $\delta_3$ & $\gamma_4$ & $\gamma_4$ \\
\hline 
\end{tabular}
\end{center}
\medskip
 
\noindent

\begin{center}
\underline{Table 1}
\end{center}
\medskip
 
\noindent

Similarly matrix units for $A_3 \cap A_3^{\ast}$  may be chosen to
standardise the map $ \psi : A_2 \cap A_2^{\ast} \rightarrow A_3 \cap A_3
^{\ast}, $ and there is a 
$\sigma \times \sigma $ matrix 
{ $
\matrrcc{a_1}{a_2}{a_4}{a_3} $}
 so that 
\medskip
 
\noindent
\[
\psi : e_{13} \otimes x \to  
  \begin{bmatrix} 0 & 0 & & & a_1 \otimes x & 0 & a_2 \otimes x
 & 0 \\ 0 & 0 & & & 0 & 0 & 0 & 0 \\
	& & 0 & 0 & a_4 \otimes x & 0 & a_3 \otimes x & 0 \\
 & & 0 & 0 & 0 & 0 & 0 & 0 \\ 
	& & & & 0 & 0 & & \\ & & & & 0 & 0 & & \\ & & & & & & 0 & 0 \\
	& & & & & & 0 & 0 \end{bmatrix}.
\]
\medskip
 
\noindent
Here $A_2$ is identified with $A(D_4) \otimes M_k$, with 
$k = \rho$, and $x \in
M_k$.  The map $\psi$ is determined by four $\sigma \times \sigma$
unitary matrices each with a $2 \times 2 $ block decomposition, 
whose entries and dimensions are indicated in Table 2.
\medskip
 
\noindent

\begin{center}
\begin{tabular}{|l|l|l|l|l|}
\hline  &&&&\\
& $u$ & $w$ & $v$ & $x$ \\
\hline  &&&&\\
$u'$ & $a_1$ & $b_1$ & $a_2$ & $b_2$ \\
\hline &&&&\\
$w'$ ~~~~~ & $d_1$ ~~~~~ & $c_1$ ~~~~~ & $d_2$ ~~~~~ &
$c_2$\\
\hline &&&&\\ 
$v'$ & $a_4$ & $b_4$ & $a_3$ & $b_3$ \\
\hline  &&&&\\
$x'$ & $d_4$ & $c_4$ & $d_3$ & $c_3$ \\
\hline  
\end{tabular}
\end{center}
\medskip
 
\noindent

\begin{center}
\underline{Table 2}
\end{center}
\medskip
 
\noindent
Thus
\[
u + v = w + x = u' + v' = w' + x' = \sigma.
\]

We have
\medskip

\[
\psi \circ \phi (e_{13}) = 
  \begin{bmatrix} 0 & 0 & & & v_1 & 0 & v_2 & 0 \\ 0 & 0 & & & 0 & 0 & 0 & 0 \\
	& & 0 & 0 & v_4 & 0 & v_3 & 0 \\ & & 0 & 0 & 0 & 0 & 0 & 0 \\ 
	& & & & 0 & 0 & & \\ & & & & 0 & 0 & & \\ & & & & & & 0 & 0 \\
	& & & & & & 0 & 0 \end{bmatrix}
\]
\medskip
 
\noindent
where
\medskip
 
\noindent
\[
\matrrcc{v_1}{v_2}{v_4}{v_3} =
 \begin{bmatrix} \alpha_1 \otimes a_1 & \beta_1 \otimes b_1& 
                 \alpha_1 \otimes a_2 & \beta_1 \otimes b_2 \\
                 \delta_1 \otimes d_1 & \gamma_1 \otimes c_1 &
                 \delta_1 \otimes d_2 & \gamma_1 \otimes c_2 \\ 
 \alpha_1 \otimes a_4 & \beta_1 \otimes b_4& 
                 \alpha_1 \otimes a_3 & \beta_1 \otimes b_3 \\
                 \delta_1 \otimes d_4 & \gamma_1 \otimes c_4 &
                 \delta_1 \otimes d_3 & \gamma_1 \otimes c_3
	\end{bmatrix}
\]
\medskip
 
\noindent
with similar expressions for the $\rho \sigma  \times \rho \sigma $ unitary
matrices
\medskip
 
\noindent
\[
\matrrcc{w_1}{w_2}{w_4}{w_3},\ \ 
 \matrrcc{x_1}{x_2}{x_4}{x_3},\ \ 
\matrrcc{y_1}{y_2}{y_4}{y_3}\ \ 
\]
\medskip
 
\noindent
arising from $\psi \circ \phi(e_{14}), \psi \circ \phi(e_{24}), \psi \circ \phi
(e_{23})$ respectively.
Thus
\medskip

\noindent
\[
\psi \circ \phi (e_{14}) = 
  \begin{bmatrix} 0 & 0 & & & 0 & w_1 & 0 & w_2 \\ 0 & 0 & & & 0 & 0 & 0 & 0 \\
	& & 0 & 0 & 0 & w_4 & 0 & w_3 \\ & & 0 & 0 & 0 & 0 & 0 & 0 \\ 
	& & & & 0 & 0 & & \\ & & & & 0 & 0 & & \\ & & & & & & 0 & 0 \\
	& & & & & & 0 & 0 \end{bmatrix}
\]
\medskip
 
\noindent
where 
\medskip
 
\noindent
\[
\matrrcc{w_1}{w_2}{w_4}{w_3} =
\begin{bmatrix} \alpha_2 \otimes a_1 & \beta_2 \otimes b_1& 
                 \alpha_2 \otimes a_2 & \beta_2 \otimes b_2 \\
                 \delta_2 \otimes d_1 & \gamma_2 \otimes c_1 &
                 \delta_2 \otimes d_2 & \gamma_2 \otimes c_2 \\ 
 \alpha_2 \otimes a_4 & \beta_2 \otimes b_4& 
                 \alpha_2 \otimes a_3 & \beta_2 \otimes b_3 \\
                 \delta_2 \otimes d_4 & \gamma_2 \otimes c_4 &
                 \delta_2 \otimes d_3 & \gamma_2 \otimes c_3	\end{bmatrix}.
\]
\medskip

\noindent {\it Proof of Lemma 3.1.}
Assume first that $\psi$ is not locally regular. (This is a consequence 
of the non local regularity of
$\eta$.)
Without loss of generality
we may assume that 
 $a_1$ is not a partial isometry. As the map $\psi \circ \phi$ is regular  
it follows that the matrix 
\medskip
 
\noindent
\[
v_1 = \matrrcc{\alpha_1 \otimes a_1}{\beta_1 \otimes b_1}{\delta_1
\otimes d_1}{\gamma_1 \otimes c_1}
\]
\medskip
 
\noindent
is a partial isometry and from this it follows that  $\alpha_1$  is
 a strict
contraction.  
Indeed if this is not the case then
the matrix unit system for $A_2$ may be chosen
so that, in addition, 
 $\alpha_1$ has the form
\medskip
 
\noindent
\[
\matrrcc{1}{0}{0}{\ast}.
\]
\medskip
 
\noindent Thus the unitary matrix
\medskip
 
\noindent
\[
\matrrcc{\alpha_1}{\beta_1}{\delta_1}{\gamma_1}
\]
\medskip
 
\noindent
has an associated form
\medskip
 
\noindent
\[
\matrrrccc{1}{0}{0}{0}{\ast}{\star}{0}{\star}{\star}.
\]
\medskip
 
\noindent
This implies that the matrix $a_1$ appears as an orthogonal part
of $v_1$, contrary to the fact that $v_1$ is a partial isometry
and  $a_1$ is not.
\medskip
 
\noindent

Similarly it follows, on consideration of the partial isometries 
\[
w_1 = \matrrcc{\alpha_2 \otimes a_1}{\beta_2 \otimes b_1}{\delta_2
\otimes d_1}{\gamma_2 \otimes c_1}
, \ \  x_1 = \matrrcc{\alpha_3 \otimes a_1}{\beta_3 \otimes b_1}{\delta_3
\otimes d_1}{\gamma_3 \otimes c_1}
, \ \ 
y_1 = \matrrcc{\alpha_4 \otimes a_1}{\beta_4 \otimes b_1}{\delta_4
\otimes d_1}{\gamma_4 \otimes c_1}
\]
that $\alpha_2$, $\alpha_3$ and $\alpha_4$ are strict
contractions. Also, by the properness of $\phi$ 
each of $\alpha_1, \alpha_2, \alpha_3, \alpha_4$ is  nonzero. 
From this it follows that none of the matrices 
of Table 1 is nonzero and none is a partial isometry.

We now show that the sixteen contractions 
$\alpha_i, \beta_i, \gamma_i, \delta_i, 1 \le i \le 4$, 
must all have the same rank.

\medskip

Since $\phi$ is star extendible we have the equations
$$\alpha_1 \alpha_1^{\ast} + \beta_1 \beta_1^{\ast} = P', \; \; \; \;
\; \alpha_1^{\ast} \alpha_1 + \delta_1^{\ast} \delta_1 = P.$$
If\ \  rank $(\alpha_1) < p' $= rank $(P' )$ \ \ then there is a rank one
projection $E \leq P'$ with $E \alpha_1 = 0$.  But then $E \beta_1
\beta_1^{\ast} E = E$, contrary to the fact that $\beta_1$ is a strict
contraction.  Thus \\  rank $(\alpha_1) = p' $.  Similarly the other
equation ensures that \ \  rank $(\alpha_1) = p =$ rank $(P)$.  For otherwise
there is a rank one projection $F \leq P$ with $\alpha_1 F = 0$ and it
follows that $\delta_1$ has an isometric part.  We conclude then that $p
= p' $.  Similar arguments, or an appeal to the symmetry of 4-cycle
algebras, leads to the equidimensionality condition
$$p = q = r = s = p{'} = q{'} = r{'} = s{'}$$
and this implies  that $\rho$ is even and the common value above is
$\rho/2$.
\medskip

Assume now that $\eta$ is not locally regular. 
Then, by the argument above,
$\psi$ satisfies the equirank condition
\medskip
$$u = v = w = x = u{'} = v{'} = w{'} = x{'} = \sigma/2.$$
We conclude then that
\medskip

\[
\mbox{ rank} (v_1) \geq \mbox{ rank } (\alpha_1 \otimes a_1) = 
\mbox{ rank } (\alpha_1) \mbox{  rank}
(a_1) = \; \displaystyle \rho \sigma /4
\]
\medskip
 
\noindent
and similarly
\medskip
 
\noindent
\[
\mbox{ rank} (v_i) \geq  \rho \sigma /4 , \ \ \ \ \ \mbox{for }1 \leq i \leq 4.
\]
\medskip
 
\noindent
Since $\psi \circ \phi$ has multiplicity $\rho \sigma $ 
it follows from the fact
that each $v_i$ is a partial isometry that 
\medskip
 
\noindent
\[
\mbox{rank} (v_i) = \displaystyle \; \rho \sigma /4 , \ \ \ \ \ 
\mbox{ for } 1 \leq i \leq 4.
\]
\medskip
 
\noindent
Thus 
$$H_1 (\psi \circ \phi) = [\mbox{rank }(v_1) - \mbox{ rank }(v_2) + \mbox{ rank
}(v_3) - \mbox{ rank }(v_4)] = 0 $$
as desired.
\hfill $\Box$
\medskip

That irregular factorisations of rigid embeddings are 
possible 
is shown in \cite{don-pow-1}. In the key stationary example there 
 $\phi$ has multiplicity 2, 
$\psi = \phi \otimes I_2$ and  $\phi$ has the matrix table
\medskip

\begin{center}
\begin{tabular}{|l|l|l|l|l|}
\hline  &&&&\\
& $1$ & $1$ & $1$ & $1$ \\
\hline  &&&&\\
$1$ & $\frac{1}{\sqrt{2}}$ & $\frac{1}{\sqrt{2}}$ & $\frac{1}{\sqrt{2}}$ & $\frac{1}{\sqrt{2}}$ \\
\hline &&&&\\
$1$ ~~~~~ & $\frac{1}{\sqrt{2}}$ ~~~~~ & $-\frac{1}{\sqrt{2}}$ ~~~~~ 
& $-\frac{1}{\sqrt{2}}$ ~~~~~ & $\frac{1}{\sqrt{2}}$ \\
\hline &&&&\\
$1$ & $\frac{1}{\sqrt{2}}$ & $-\frac{1}{\sqrt{2}}$ & $-\frac{1}{\sqrt{2}}$ & $\frac{1}{\sqrt{2}}$ \\
\hline  &&&&\\
$1$ & $\frac{1}{\sqrt{2}}$ & $\frac{1}{\sqrt{2}}$ & $\frac{1}{\sqrt{2}}$ & $\frac{1}{\sqrt{2}}$ \\
\hline
\end{tabular}
\end{center}
\medskip

\begin{lma}
 Let $\phi : A_1 \rightarrow
A_2,$ $\psi: A_2 \rightarrow A_3$  be locally regular
star extendible embeddings between 4-cycle algebras.
If the composition $\psi \circ \phi$ is a rigid embedding
then $\phi, \psi$ are rigid embeddings.
\end{lma}
\medskip 

\begin{pf}
We may assume that $A_1 = A(D_4)$ and that the maps are unital.
The general case follows readily from this one.
By hypothesis the matrices of Tables 1 and 2 are partial isometries.
Since $\psi \circ \phi$  is rigid it follows that $v_1$ and $w_2$ have the
same final projections and  hence, in particular, the matrices
\[
[\alpha_1 \otimes a_1 \ \  \beta_1 \otimes b_1] \ \ \ \mbox{and} \ \ \  
[\alpha_2 \otimes a_2 \ \  \beta_2 \otimes b_2 ]
\]
have the same final projections. The partial isometries $a_1$ and $a_2$ 
have orthogonal final projections  
since they appear in the unitary matrix 
$\bigl[ \begin{smallmatrix} a_1 & a_2 \\ a_3 & a_4 \end{smallmatrix} \bigr]$
and so
it follows that the pair
$\alpha_1 \otimes a_1$ and
$\beta_2 \otimes b_2 $
 and hence the pair $a_1, b_2$,  have coincidental final projections.

Similarly
$b_2, c_3$ have the same initial projection, 
$c_3, d_4$ have the same final projection and
$d_4, a_1$ have the same initial projection. By the star extendibility of $\psi$ 
we have the equality $d_4^*c_1b_2^*a_1 = a_1^*a_1$ and
it follows that the matrices $a_1, b_2, c_3, d_4$ 
determine a direct summand of $\psi$ which is a rigid embedding.
In fact this summand is (inner) unitarily equivalent to a direct sum 
of $k$ equivalent 
rank one rigid embeddings where $k = $rank$(a_1)$.

Continuing in this way it follows that $\psi$ is a rigid embedding with $H_1(\psi)
= [\delta ]$
where
\[
\delta = \mbox{rank}(a_1) - \mbox{rank}(b_1) + \mbox{rank}(c_1) 
- \mbox{rank}(d_1).
\]
Since $\psi $ and $\psi \circ \phi$ are rigid embeddings
it follows that
$\phi$ is also a rigid embedding.
\end{pf}
\medskip

\begin{thm}
For generic
limit algebras of 4-cycle algebras
the limit homology group is a well-defined invariant
for star extendible isomorphism.
\end{thm}
\medskip

\begin{pf}
Suppose that $\{A_k, \alpha_k\}, $ $\{A_k', \alpha_k'\}$
are generic systems with limits $A$ and $A'$ and that $A$ and $A'$
are star extendibly isomorphic.
Then there is a star extendible commuting diagram isomorphism between
$\{A_k, \alpha_k\} $ and  $\{A_k', \alpha_k'\}$. 
Composing embeddings we may assume
that the crossover maps are proper. This is so because the systems are generic.
By Lemmas  3.1 and 3.2 either $\{A_k, \alpha_k\}$ and $\{A_k', \alpha_k'\}$
are regularly isomorphic or both $H_1(\{A_k, \alpha_k\})$
and $H_1(\{A_k', \alpha_k'\})$ are zero.
\end{pf}
\medskip

\begin{remark} A similar result holds for limits of $2n$-cycle algebras
for $n \ge 3$ details of which will appear in \cite{don-pow-3}. 
It is curious that the complexity of the embeddings between 
such higher order cycle algebras is
offset by the fact that they can be shown to be automatically
locally regular. Thus one can then move to the consideration
of a couterpart to Lemma 3.3 and bypass entirely the detailed analysis
of irregular factorisations of regular embeddings that we have given here
and in \cite{don-pow-1}. 
\end{remark}

\section{$H_1$ classifications.}

If $\alpha$ is a star extendible embedding between
connected digraph algebras then its multiplicity is defined
to be the multiplicity of its star algebra extension.
Alternatively, the multiplicity of $\alpha $ in the case of a  regular
embedding is $|r|$ where $H_0(\alpha)$ is
realised as multiplication by the
integer $r$.

A  direct system of  digraph algebras,
with star extendible embeddings, 
is said to be an  {\it odd}  system if only finitely many
of its embeddings have even multiplicity.
\medskip

\begin{thm}
A  star extendible 
isomorphism between 
{\it odd} generic direct systems of 4-cycle algebras
or between systems with nonzero $H_1$ invariant
is necessarily a regular isomorphism.
\end{thm}

\begin{pf} 
The homology of 
the maps of an odd generic system
are eventually nonzero
and so the argument of Theorem 3.3 applies.
\end{pf}
\medskip

Star extendible isomorphisms
between the
 algebraic direct limits are automatically regular isomorphisms
and so  the algebraic limit algebras possess all the invariants that
their systems possess for regular isomorphism. 

It is still the case that generic even systems with trivial
$H_1$ invariant are regularly isomorphic if 
they are star extendibly isomorphic,
but in the absence of automatic regularity this
regular isomorphism must be constructed.
\medskip
 
\begin{thm}  Let $\ca, \ca{'}$ be generic systems of 4-cycle
algebras with algebraic direct limits
$A = \algindlimit \ca \ , \  A{'} =
\algindlimit \ca{'}$
respectively.  Then 
$A$ and $A{'}$ are star extendibly
isomorphic if and only if
$\ca$ and $\ca{'}$ are regularly
isomorphic.
\end{thm}
\medskip

\begin{pf}
By the analysis of the last section 
we may assume that $\ca, \ca{'}$ have embeddings $\alpha_k, \alpha_k'$
which satisfy the equidimensionality conditions. In particular 
$\ca, \ca{'}$ are of matroid type and have trivial $H_1$ invariant.
Thus $\ca = \{A_k, \alpha_k\}, \ca' = \{A_k', \alpha_k'\}$
and $K_0\alpha_k , K_0\alpha_k'$ have the form
\[
 \left[ \begin{array}{cccc}
t_{k}&t_k&0&0\\
t_k&t_k&0&0\\
0&0&t_k&t_k\\
0&0&t_k&t_k
\end{array} \right], \ \ \ 
 \left[ \begin{array}{cccc}
u_{k}&u_k&0&0\\
u_k&u_k&0&0\\
0&0&u_k&u_k\\
0&0&u_k&u_k
\end{array} \right].
\]
Since $A$ and $ A'$ are star extendibly isomorphic it follows that
there is a scaled group isomorphism
between the direct systems $K_0\ca, K_0\ca'$
in the form of a commuting diagram.
Moreover we may assume that
he crossover maps have the form
\[
 \left[ \begin{array}{cccc}
v_{k}&v_k&0&0\\  
v_k&v_k&0&0\\
0&0&v_k&v_k\\
0&0&v_k&v_k
\end{array} \right]
\]
where $v_k$ is even for all $k$. 
(If this is not already so then the crossover maps have the form
\[
 \left[ \begin{array}{cccc}
a&b'&0&0\\
b&a'&0&0\\  
0&0&c&d'\\  
0&0&d&c'    
\end{array} \right]
\]
where, by star extendibility,  $a+b = a'+b' = c+d = c'+d'$. Composing
such a map with a given $K_0$ map gives a new crossover map
of the desired type.)
Thus each crossover 
map has a lifting to a generic
embedding with zero $H_1$ map. 
Construct liftings of the crossover maps, in order, to rigid embeddings
with trivial homology,
and use $K_0H_1$
uniqueness to arrange commuting triangles.
\end{pf}
\medskip

\noindent {\bf Remarks.} 1. \ In fact the generic requirement can be dropped
and Theorem 4.2 holds for general rigid systems of 4-cycle algebras.
As Allan Donsig has observed, it is straightforward to obtain
the automatic regularity of the factors
$\phi, \psi$ of a composition $\psi \circ \phi$ which is of nongeneric type.
This follows essentially from the fact that a partial isometry with
support in three block subspaces is necessarily regular.
 
2.   The theorem above implies a simplification of the isomorphism
problem for algebraic limits in the following sense.

Choose partial matrix unit systems $\{e_{ij}^k\}$ for $A_k$ for 
$k = 1, 2, \dots$ in the usual way so 
that each  $e_{ij}^k$ is a sum of some of the 
matrix units of the system $\{e_{ij}^{k+1}\}$. Whilst the semigroup
$ S = \{e_{ij}^k\}$ formed by the totality of all these matrix units depends
on the system $\ca$, the semigroup ring
$R_A = {\Bbb Z}[S]$ is, by the theorem, a well-defined invariant for
$A$. Thus the algebra $A$ comes with a canonical ring inclusion
$
R_A \to A $  for which $ A = {\Bbb C} \otimes_{\Bbb Z} R_A
$.

3. It is likely that standard techniques 
lead to the fact that $A$ and $A'$ are star extendibly isomorphic
if and only if they are isomorphic as complex algebras.
Alternatively, although more indirectly,
in view of the classifications below would be enough
to formulate the $K_0H_1$ invariants as algebra isomorphism 
invariants.
\medskip
 
4. It is still an open problem whether the conclusion of 
Theorem 4.2 holds for the
limit algebras of arbitrary regular star extendible 
systems of digraph algebras.
\medskip

\begin{defn} 
For a generic index 4 matroid algebra position
$\pi : D_1 \oplus D_2 \to B$, or, more generally, for a generic index 4
position of 2-symmetric algebras,  the partial isometry homology group 
$H_1(\pi)$ is defined to be the abelian 
group $H_1(A)$ where $A$ is the associated 4-cycle
limit algebra. 
\end{defn}
\medskip

Let  $\zeta  : F \oplus F \to F$,
with $F = \algindlimit M_{2^k}$, 
be the unital 
inclusion determined by the stationary
system given in Section 2, based on the simplest generic embeddings.
It follows from the argument of the last proof
 that $\zeta$ is the unique position
with trivial homology invariant.
 \medskip

\begin{thm}
Up to position
there is a unique generic index 4 matroid algebra inclusion
\[
\zeta  : F \oplus F \  \to \ \ F
\]
\noindent
with $H_1(\zeta) = 0.$
\end{thm}
\medskip

This uniqueness
might suggest that 
$H_1(\pi)$
provides a complete invariant for the generic index 4 inclusions
but this is not the case.
The most transparent obstacle
is given by a {\it coupling class} in ${\Bbb R}/{\Bbb Q}$. 
If $A = \algindlimit {(A_k ,\alpha_k)}$ 
then this class is defined to be that which is determined
by any of the products
\[
\kappa{(A)} = \prod_{k=l}^{\infty} \frac{|\delta_k |}{p_k},
\]
for suitably large $l$, where $p_k$ is the multiplicity 
of $\alpha_k$ and $H_1(\alpha_k) = [\delta_k]$.
\medskip       

The following partial classification 
also follows from the more general Theorem 5.3.
\medskip
 
\begin{thm}
Let $\pi_1, \pi_2$ be generic index 4 matroid algebra inclusions
of $D_1 \oplus D_2$ in $D$ and assume that $\kappa(\pi_1) = \kappa(\pi_2) = 0$.
Then $\pi_1$ and $\pi_2$ are conjugate if
and only if the abelian groups $H_1(\pi_1)$
and $H_1(\pi_2)$ are isomorphic.
\end{thm}
\medskip       

\begin{pf}
The necessity of the condition for conjugacy follows from
Theorem 3.3. Assume then that  $H_1(\pi_1)$
and $H_1(\pi_2)$ are isomorphic.
As in the proof of Theorem 4.2 there is a commuting diagram isomorphism
between the direct systems $K_0\ca$ and $K_0\ca '$ associated with
$\pi_1$ and $\pi_2$.
Lift the first crossover map of this isomorphism to a rigid embedding
$\phi : A_1 \to A_k'$ for some $k$.
Under the hypotheses it is possible to lift
the next crossover map, $Y_1 : K_0A_k' \to K_0A_l$ say, after
increasing $l$ if necessary, to a rigid embedding $\psi_1$ so that
$\psi_1 \circ \phi_1$ is equal to the given embedding
$i : A_1 \to A_l$. This can be seen from  the essential part 
of the proof of Theorem 11.22 in
\cite{pow-book}. Thus, by $K_0H_1$ uniqueness it is enough to choose $l$ 
large enough and
$\delta$ in the homology range of $Y_1$ so that
$\delta H(\phi_1) = H_1(i)$. This is possible since, for 
large $l$ the ratios 
$|H_1(i)|/|H_0(i)|$ are arbitrarily small.

One can continue in this way to obtain a commuting diagram of rigid embeddings
which determines the desired conjugacy.
\end{pf}

\section{$K_0H_1$ classifications}

The invariants $H_1(\pi)$ and $\kappa(\pi)$ take 
no account of the order of the
summands and so, as formulated, they cannot serve as invariants
for the conjugacy of
automorphisms with summand respecting restrictions.
To determine such conjugacy we now consider $K_0H_1$ invariants for the ordered
inclusion and, most decisively, the joint scale in $K_0H_1$
introduced in 
 Donsig and Power
\cite{don-pow-2}.

Let $\{A_k, \alpha_k\}$ 
be a generic 4-cycle algebra system with algebraic direct limit
$A$. The identification of the 
scaled $K_0$ group invariant $(K_0A, \Sigma_0A)$
and the associated classification of the algebras $A\cap A^*
= D_1 \oplus D_2$ has
been considered in detail in 
\cite{don-pow-2}.
If  $\{A_k, \alpha_k\}$ is of matroid type then $D_1$ and $D_2$
are stably isomorphic even matroid algebras and are isomorphic
in the unital
equidimensional case.
Thus $K_0(A)$ may be realised naturally as a 
subgroup of ${\Bbb Q} \oplus {\Bbb Q}$, with the usual product
order and with scale determined by a product of 
(possibly infinite) intervals,
$I_1 \times I_2$.
The algebra $A$ is unital if and only if both $I_1$ and $I_2$ are finite closed
intervals. See Dixmier \cite{dix}.

More generally $K_0\alpha_k$ has the form
\[
 \left[ \begin{array}{cccc}
a_{k}&b_k&0&0\\
b_k&a_k&0&0\\
0&0&c_k&d_k\\
0&0&d_k&a_k
\end{array} \right]
\]
where $ a_k + b_k = c_k + d_k$. Let $p_k =  a_k + b_k$ let \ 
$q_k = a_k - b_k$ (and
without loss of generality assume that $q_k \ge 0 $ for all $k$) and
let $G_p$ and $G_q$ be the subgroups of ${\Bbb Q}$ associated
with the generalised integers
$p = p_1p_2\dots,\  q = q_1q_1\dots$.
One can show that $K_0D_1$ can be realised
as $G_p \oplus G_q$ in this case. In the odd case $K_0D_1$ is generated by 
$G_p \oplus G_q$ and $(1/2, 1/2)$.

Note that the scale
$\Sigma_0A$ can be identified as the set of
$K_0$ classes of projections $\phi(1)$
associated with the star extendible injections
\[
\phi : {\Bbb C} \to A_k.
\]
Likewise 
the scale $\Sigma_1A$ of $H_1A$ can be defined as the set of elements
$(H_1\psi)(g)$ associated with the morphisms
\[
H_1\psi : H_1(A(D_4)) \to H_1A_k
\]
induced by regular embeddings
\[
\psi  : A(D_4) \to A_k
\]
for some $k$, where $g$ is a fixed generator for $H_1(A(D_4))$.

The joint scale admits a similar definition.
\medskip

\begin{defn}
Let $A$ be a generic limit of 4-cycle algebras. Then 
the joint scale $\Sigma A$ of $ K_0A \oplus H_1A$ is defined to be the subset 
of $\Sigma_0A \times \Sigma_1A$ consisting of elements
\[
(K_0\phi([e_{11} \oplus e_{33}]), (H_1\phi)(g))
\]
associated with the  rigid  embeddings
\[
\phi : A(D_4) \to A_k
\]
for some $k \ge 1$.
\end{defn}

The next theorem, the lifting Theorem 5.2, 
can also be obtained from Theorem 4.2 above and
the main result in Donsig and Power \cite{don-pow-2}. For completenes
we present a proof.

In \cite{don-pow-2} we have discussed the circumstances under which
it is possible to lift a scaled group homomorphism
\[
\gamma_0 \oplus \gamma_1 : K_0A_1 \oplus H_1A_1 \to K_0A_2 \oplus H_1A_2
\]
to a rigid embedding between the 4-cycle algebras $A_1, A_2$.
(We assume that $\gamma_0$ respects the ordered summand structure of 
$K_0$, so that $\gamma_0$ is necessarily block diagonal.)
A necessary and sufficient condition for a lifting to exist
is that $\gamma_0$ be of rigid type, that is, implemented by a 
matrix of the form
\[
\left[ \begin{array}{cccc}
a&b&0&0\\
b&a&0&0\\
0&0&c&d\\
0&0&d&c
\end{array}\right]
\]
with $a+b = c + d$, and that $\gamma_0 \oplus \gamma_1$
maps $\Sigma A_1$ into $\Sigma A_2$.
The rigid type nature of $\gamma_0$ is also equivalent to symmetry preservation,
by which we mean that $\gamma_0 \circ \theta = \theta \circ \gamma_0$, 
for each of the (four) canonical symmetries $\theta$
of the ordered 
$K_0$ group.

In general symmetry preservation is not an automatic
consequence of joint scale preservation. This can be seen for example
for the scaled group homomorphism
\[
\gamma_0 \oplus \gamma_1 = 
\left[ \begin{array}{cccc}
6&6&0&0\\
2&2&0&0\\
0&0&6&6\\
0&0&2&2  
\end{array}\right] \oplus 0
\]
viewed as a map from $K_0(A(D_4))$ to $K_0(A(D_4) \otimes M_{12})$.
However, one can readily check that 
joint scale preservation does ensure that the column sums 
of $\gamma_0$ are constant (and hence that $\gamma_0$ extends to 
a scaled group homomorphism from $K_0C^*(A_1)$ to $K_0C^*(A_2)$).
Thus, if in fact $\gamma_0$ is a map from
$K_0(A(D_4) \otimes M_r)$ to $K_0(A(D_4) \otimes M_s)$ which 
preserves the order unit (and so maps $(r,r,r,r)$ to $(s,s,s,s)$) then the
row sums as well as the column sums coincide, and so $\gamma_0$
is automatically of rigid type in this case.

Naturally we say that a (summand respecting scaled group)
isomorphism $\gamma_0 : K_0A \to K_0A'$ is of rigid type if it pulls back
to a commuting diagram of rigid type embeddings.
This condition is similarly
equivalent to symmetry preservation for the canonical
symmetries on $K_0A, K_0A'$. (See also the discussion in \cite{pow-book}.)

These remarks explain why one can drop the symmetry 
preservation hypothesis in the next theorem in the case when
the algebras are unital and are determined by equidimensional
systems.

\begin{thm}
Let $A, A'$ be  (algebraic) limits of 4-cycle algebras with respect to
rigid embeddings and let
\[
\gamma_0 \oplus \gamma_1 : K_0A \oplus H_1A \to K_0A' \oplus H_1A'
\]
be an abelian group isomorphism where $\gamma_0$ 
is a scaled group isomorphism, which in the nonunital case
is symmetry preserving.
Then the following assertions are equivalent.

(i) There is a star extendible isomorphism $\phi : A \to A'$
for which $\gamma_0 \oplus \gamma_1 = K_0\phi \oplus H_1\phi$.

(ii)  
The isomorphism 
$\gamma_0 \oplus \gamma_1 $
effects a bijection between the joint scales.
\end{thm}

\begin{pf} If (i) holds then by Theorem 4.2 and the remarks
that follow it there is a regular
star extendible isomorphism
between the systems $\{A_k, \alpha_k\}$ and
$\{A_k', \alpha_k'\}$  for $A, A'$ with the same induced $K_0$ map as $\phi$.
Such an isomorphism induces a $K_0 \oplus H_1$ map 
$\gamma_0 \oplus \gamma_1$, with $\gamma_0$ symmetry preserving,
which effects a bijection
between
the joint scales. 

Assume now that (ii) holds.
Without loss of generality assume that $A_1 = A_1' = A(D_4)$
and consider a rigid embedding 
$\psi : A(D_4) \to A_1$ which determines
the triple
\[
([p] \oplus [q], \delta)
\]
where $p = \psi(e_{11}), q = \psi(e_{33})$ and $ \delta = (H_1\psi)(g)$.
By the joint scale preservation there is a
generic embedding $\eta : A(D_4) \to A_k'$, for some $k'$, such
that
\[
\gamma_0 \oplus \gamma_1(  ([p]\oplus[q],\delta)  = 
([\eta(e_{11} \oplus e_{33})],(H_1\eta)(g)).
\]
Since $\gamma_0$ is symmetry preserving, $\gamma_0$ and $K_0\eta$
agree and $\eta$ is a lifting. (See also Lemma 11.4 of \cite{don-pow-2}.)

Consider now a natural copy
of $A(D_4)$ in $A_k'$, that is, any copy for which the inclusion
is a multiplicity one rigid embedding. The associated injection
$\eta : A(D_4) \to A_k'$ determines an element
$\gamma_0^{-1} \oplus \gamma_1^{-1}(([\eta(e_{11} \oplus e_{33})], H_1\eta(g)))$
in the joint scale of $K_0A \oplus H_1A$ and so, for some generic embedding
$\xi : \eta(A(D_4)) \to A_l$ we have
\[
\gamma_0^{-1} \oplus \gamma_1^{-1}([\eta(e_{11} \oplus e_{33})],(H_1\eta)(g))
 = ([\xi(e_{11} \oplus e_{33})], H_1\xi(g)).
\]
By symmetry preservation $K_0\xi$ agrees with the restriction
of $\gamma_0^{-1}$.

We claim that the embedding $\xi$ has an extension
$\xi : A_k' \to A_l$, perhaps after increasing $l$, which is also generic.
Indeed, since $\gamma_0^{-1}$ is a scaled group isomorphism
we can first extend the restriction
$\xi|(\eta(A(D_4)) \cap (\eta(A(D_4)))^*)$
to a C*-algebra injection
\[
\xi_* : A_k' \cap (A_k')^* \to A_l \cap (A_l)^*
\]
so that $K_0\xi_*$  agrees with 
$\gamma_0^{-1}|(A_k' \cap (A_k')^*)$.
Each matrix unit $e$ in $A_k'$ (for some fixed matrix unit choice)
admits a unique factorisation
$e = e_1fe_2$ with $f$ a matrix unit in $\eta(A(D_4)) $
and $e_1, e_2$ matrix units
in $A_k' \cap (A_k')^*$. Thus $\xi$ and $\xi_*$ have a joint
extension $\xi : A_k' \to A_l$ and this extension is also generic.

Now $ \xi \circ \eta : A_1 \to A_l$ has the 
same $K_0 \oplus H_1$ map as the given embedding
$\alpha_{l-1} \circ \alpha_{l-2} \circ \dots \circ \alpha_1$
and so by $K_0H_1$ uniqueness we may replace $\xi$ by an inner
conjugate embedding to get equality.
Continue this process to get a regular isomorphism
between the systems and hence a star extendible isomorphism
between $A$ and $A'$.
\end{pf}

Now let $\pi : D_1 \oplus D_2 \to D$ be the ordered subalgebra position
determined by the algebra $A$ as in the last theorem.
Define $K_0H_1(\pi)$
to be the group $K_0A \oplus H_1A$ and let
$\Sigma(\pi)$
be the joint scale $\Sigma A$ of $K_0A \oplus H_1A$. Since $A$ is determined by the
ordered inclusion $\pi$ the pair $(K_0H_1(\pi), \Sigma(\pi))$
is an invariant for ordered inclusion. The next theorem follows immediately
from Theorem 5.2 and the remarks above and resolves the problem of classifying
generic inclusions (of unspecified index).

\begin{thm}
Let $D_1, D_2$ be 2-symmetric algebras and let $\pi_1, \pi_2$ be 
generic ordered inclusions of 
$D_1 \oplus D_2$ in $D$. Then there is a star isomorphism $\alpha : D \to D$
with 
$\alpha (\pi_1(D_i)) = \pi_2(D_i)$ for $ i = 1,2$
if and only if 
there is a group isomorphism 
\[
\gamma_0 \oplus \gamma_1 : K_0H_1(\pi_1) \to K_0H_1(\pi_2),
\]
where $\gamma_0$ is a symmetry preserving scaled group isomorphism 
and
\[
 \gamma_0 \oplus \gamma_1(\Sigma(\pi_1)) = \Sigma(\pi_2).
\]
\end{thm}

It is also of interest to determine the symmetries and automorphisms
of these positions and for this one must consider realisations 
of the invariants.
We give an indication of this determination in the remainder of this section.

The appropriate outer automorphism group is given in the following 
definition.
The term approximately inner automorphism
 indicates automorphisms which are pointwise limits of a sequence of
inner unitary automorphisms. 
\medskip

\begin{defn}
Let 
$\pi : D_1 \oplus D_2 \to D$
be a 2-symmetric algebra inclusion associated with a generic 4-cycle
algebra system. Then $Out_{D_1 \oplus D_2}(D)$ is 
the group of automorphisms of $D$ which leave invariant the subalgebra
$D_1 \oplus D_2$ and for which the restriction to $D_1 \oplus D_2$ 
is an automorphism which is approximately inner.
\end{defn}

\begin{thm}
If $\pi$ is as above and 
 $H_1(\pi) \ne 0$ then $Out_{D_1 \oplus D_2}(D)$ is isomorphic
to the group of symmetry preserving automorphisms 
$\gamma_0 \oplus \gamma_1$ of 
$(K_0H_1(\pi), \Sigma(\pi))$.
\end{thm}

\begin{pf}
In view of Theorem 5.3 it is enough to show that a subalgebra 
respecting automorphism
which induces the identity map on $(K_0H_1(\pi), \Sigma(\pi))$
has a restriction which is approximately inner. However 
by $K_0H_1$ uniqueness
this is routine.
\end{pf}
\medskip

\begin{thm}
Let $D_1$ and $D_2 $ be even unital matroid algebras and let $\pi :
 D_1 \oplus D_2 \to D$
be a generic index 4 inclusion with $H_1(\pi) \ne 0$.. Then 
$Out_{D_1 \oplus D_2}(D)$
is equal to $Aut(H_1(\pi ))$ if $\kappa (\pi ) = 0$ and is equal to ${\Bbb Z}_2$
otherwise.
\end{thm}
\medskip

\begin{pf}
Let $A$ be the usual 4-cycle algebra limit algebra for which $A \cap A^* = 
 D_1 \oplus D_2$.
It follows from Theorem 7.2 of \cite{don-pow-2}
that $\Sigma A$ is invariant under the homology inversion
$id \oplus -1$
and that $\Sigma A$ splits as a direct sum, say 
$\Sigma_{00}A \oplus \Sigma_1 A$.
If $\kappa (\pi ) > 0$ then $\Sigma_1A$ is a finite interval in $H_1A $
and if $\kappa (\pi ) = 0$ then $\Sigma_1A = H_1A$. Now Theorem 5.5
completes the proof, since symmetry preservation is automatic in the unital
case.
\end{pf}

It is intriguing that in the 2-symmetric case there may exist obstructions to 
homology inversion as we now show.

Let $\alpha_1 : A_1 \to A_2$ be a rigid embedding between 4-cycle
algebras. Then the {\it abstract homology range} of $\alpha_1$, 
denoted $hr(\alpha_1) $
is the set of maps $H_1(\beta)$ for which $\beta : A_1 \to A_2$  is rigid
with $K_0(\beta) = K_0(\alpha_1)$.
Thus if , for example, $\alpha_k$ has $K_0$ map
\[
 \left[ \begin{array}{cccc}
n_{k}&2&0&0\\
2&n_k&0&0\\
0&0&n_k&2\\
0&0&2&n_k
\end{array} \right]
\]
then $\alpha_k $ has possible signatures 
$\{n_k,0,2,0\}$,
$\{n_k-1,1,1,1\}$,
$\{n_k-2,2,0,2\}$
with respective homology maps
$[n_k+2], [n_k-2] , [n_k-6]$. Thus the homology range is 
$\{n_k+2, n_k-2, n_k-6\}$.
\medskip

\begin{defn}
\cite{don-pow-2} \ \ The generic 4-cycle algebra system
$\{A_k, \alpha_k\}$ is hr-symmetric
if for each pair $j, i $  the homology
range of $ \alpha_j \circ \dots  \circ \alpha_i$ lies in ${\Bbb Z}_+ $ or 
${\Bbb Z}_-$ and does not contain $0$.
\end{defn}
\medskip

Such a system may be constructed as follows.
Returning to our earlier notation for
$\alpha_k$ note that $K_0(\alpha_k)$ is defined by
the triple $(p_k,q_k,r_k)$ where
$p_k = a_k + b_k, q_k = a_k - b_k,
 r_k = c_k - d_k$. (Again, we assume that $c_k \ge d_k$ for all $k$.)
 The triple for the composition
 $\alpha_{k+1} \circ \alpha_{k} \circ \dots \circ \alpha_1$
 is the triple
 \[
 (P_k,Q_k,R_k) = (p_kp_{k-1} \dots p_1, q_kq_{k-1} \dots q_1,r_kr_{k-1}
  \dots r_1)
 \]
 and the homology range of the composition is contained in the interval
 $[s,t]$ where
 \[
 t = P_k - |Q_k - R_k|. \ \  \ 
 s = Q_k + R_k - P_k.
 \]
 In particular if we choose $\alpha_k$ as above with
 $n_k$ increasing, such that for all $k$

 \[
\frac{2(n_k-2)}{n_k+2} > \prod_{i=1}^{k-1} \frac{n_i+2}{n_i -2},
\]
then the homology range of 
$\alpha_{k+1} \circ \alpha_{k} \circ \dots \circ \alpha_1$
is positive for all $k$.

\begin{prop}
Let $A$ be the generic limit algebra of an hr-asymmetric
system. Then the joint scale $\Sigma A$ in $K_0 \oplus H_1$ is
not invariant under the homology inversion
$\gamma_0 \oplus \gamma_1 = id \oplus -1$.
\end{prop}

\begin{pf}
Adding a rank one rigid embedding $A(D_4) \to A_1$ we may asssume that
$A_1 = A(D_4)$.
Thus the element $[e_{11} \oplus e_{33}] \oplus g$ in $\Sigma A_1$
determines an element of the joint scale of $A$.
Suppose that 
$[e_{11} \oplus e_{33}] \oplus -g$
 is also an element of the joint scale,
appearing as an element of $K_0A_k \oplus H_1A_k$
for some $k$.
Then we conclude that both $g$ and $-g$ belong to the homology range of
$\alpha_{k-1} \circ \dots \circ \alpha_1$ contrary to hypothesis.
\end{pf}

\medskip

\begin{remark} Another possible obstacle 
to homology inversion, in the odd 2-symmetric case,
is a mod 4 congruence class which may be implicit in the joint scale.
This is considered in detail in \cite{don-pow-2}. This 
somewhat subtle obstruction 
is purely arithmetic and 
in contrast to hr-asymmetry it is annihilated
by tensoring with an even
matroid algebra.
\end{remark}

\section{Operator Algebras}

We now obtain approximate versions of the lemmas of section 3
and we retain the notation from that section.
The following definition will be useful in arguments involving
embeddings which are almost locally regular.
\medskip

\begin{defn}  A star extendible homomorphism $\phi$ 
between 4-cycle algebras is $\epsilon$-strict if
each matrix entry $\alpha_i, \beta_i, \gamma_i, \delta_i$, for $1 \leq
i \leq 4$, is distance at most $\epsilon$ from a partial isometry. 
\end{defn}
\medskip

Let $\phi: A(D_4) \rightarrow A_2$ be a general star extendible
embedding and let $\phi' = \alpha  \circ \phi$ where $\alpha :
A_2 \rightarrow A_3$ is a generic embedding.  Let  $\alpha_i',  \beta'_i
\dots$, be the matrices corresponding to $\phi'$ 
as in Table 1.  
Since $\alpha$ is generic it follows that for each $1 \leq i \leq 4$,
$$|| \alpha'_i|| = || \beta'_i|| = || \gamma_i'|| = || \delta_i'|| =
\mbox{max} \left\{ || \alpha_i ||, || \beta_i ||, || \gamma_i ||, ||
\delta_i || \right\}$$
A 4-cycle algebra embedding with this property is said to be
 {\it norm-symmetric}.
\medskip

\begin{lma}
 Let $\left[ \begin{array}{cc} x & \epsilon_1 \\ \epsilon_2 &
y \end{array} \right]$ be a contraction which is $\eta$-close to a partial
isometry.  If $\delta =  \eta +||\epsilon_1||+||\epsilon_2||$, and 
$\delta < \frac{1}{8}$, then
$x$ is $8\delta$-close to a partial isometry.
\end{lma}
\medskip  

\begin{pf}
The matrix
\[
 \matrrcc{x^{\ast}x + \varepsilon_2^{\ast} \varepsilon_2}
{x^{\ast}\varepsilon_1 + \varepsilon_2^{\ast}y}
{\varepsilon_1^{\ast}x+y^{\ast}\varepsilon_2}
{\varepsilon_1^{\ast} \varepsilon_1 + y^{\ast} y}
\medskip
 \]
is $ 2 \eta$-close to a projection and so $ \matrrcc{x^{\ast}
x}{0}{0}{y^{\ast}y}$  is $ 2 \delta$-close  to a projection.
In particular  $
||(x^{\ast}x)^2 - (x^{\ast}x)|| \le 4\delta $ from which it follows,
if $\delta \leq \frac{1}{2} $, that 
$ x $ is
$ 8 \delta$-close to a partial isometry.
 \end{pf}
\medskip

\begin{lma}
 Let $\epsilon < \frac{1}{16}$.  Let $\phi : A (D_4) \rightarrow A_2,
\; \; \psi : A_2 \rightarrow A_3$ be proper embeddings of 4-cycle algebras for
which the composition $\psi \circ \phi$ is $\epsilon^2$-strict.  
If $\psi$ is not
$\epsilon$-strict and $\phi$ is norm-symmetric then $\phi$ is not
$(\epsilon/51)^2$-strict.
\end{lma}
\medskip

\begin{pf}
Since $\psi$ is not $\epsilon$-strict at least one of the matrices of
Table 2 is not $\epsilon$-close to a partial isometry.  Without loss of
generality we may assume this matrix to be $a_1$.  By assumption the
matrix
\[
v_1 = \matrrcc{\alpha_1 \otimes a_1}{\beta_1 \otimes b_1}{\delta_1  
\otimes d_1}{\gamma_1 \otimes c_1}
\]
\medskip

\noindent is $\epsilon^2$-close to a 
regular partial isometry.  We now deduce from
these two facts that $\alpha_1$ has norm no greater than $1 -
(\epsilon /50)^2$.

Let $t = || \alpha_1 ||$.  This norm is attained and there is a
block decomposition
$$\alpha_1 = \left[ 
\begin{array}{cc} t & \epsilon_1 \\ \epsilon_3 & \star \end{array}
\right]$$
and an associated induced decomposition of $v_1$;
$$v_1 = \left[ \begin{array}{ccc} t \otimes a_1 & 
\epsilon_1 \otimes a_2 & \epsilon_2
\otimes b_1 \\ \epsilon_3 \otimes a_1 & \star & \ast \\
 \epsilon_4 \otimes d_1 & \ast &
\ast \end{array} \right].$$
Since $\psi \circ \phi$ is $\epsilon^2$-strict it follows that
$$
\left[ \begin{array}{cc} t \otimes a_1 & \epsilon_1 \otimes a_1 \\ \epsilon_3
\otimes a_1 & \ast \end{array} \right]
$$
is $\epsilon^2$-close to a partial isometry and so, in view of the last lemma,
the matrix $ta_1 = t \otimes a_1$ is $8(\epsilon^2 + ||\epsilon_1 || + 
||\epsilon_2||)$
close to a partial isometry.  We have 
$$
\begin{array}{lll}
8 (\epsilon^2 + ||\epsilon_1|| + ||\epsilon_2||) & \leq & 
8 \; (\epsilon^2+2 \; (1-t^2)^{\frac{1}{2}})\\
& \leq & 8 \epsilon^2 + 16 \sqrt{2} \; (1-t)^{\frac{1}{2}}. \end{array}
$$
Since $a_1$ is $(1-t)$-close to $ta_1$, and yet not $\epsilon$-close to a
partial isometry, it follows that
$$\epsilon \leq (1-t)+ 8\epsilon^2 + 16 \sqrt{2} \; (1-t)^{\frac{1}{2}}.$$
Thus
\[
\epsilon(1- 8 \epsilon)  \leq  (1+16 \sqrt{2}) (1-t)^{\frac{1}{2}}
\]
and so
\[
t  \leq 1 - \epsilon^2(2(1+16 \sqrt{2}))^{-2}\\
 \leq 1 - (\epsilon /50)^2.
\]
\medskip

Considering $w_1, x_1$ and $y_1$ in similar ways it follows that $||
\alpha_i|| \leq 1 - (\epsilon/50)^2$ for $i=2,3,4$.  Since $\phi$ is
norm-symmetric it follows that $\phi$ cannot be 
$(\epsilon/51)^2$-strict.  
Indeed this would imply that $\alpha_1$ and $\beta_1$ are
$(\epsilon/51)^2$-close to partial isometries, one of which at least is
necessarily nonzero, and thus of norm 1.
\end{pf}
\medskip

The next lemma is a partial generalisation of the factorisation
dichotomy of Lemma 3.1.
\medskip

\begin{lma}
 Let $\epsilon < \frac{1}{16}$ and let $\epsilon_{\ast} \leq$
$\frac{1}{17}(\epsilon^4/51^6)$.  
Let $\phi : A_1 \rightarrow A_2, \psi: A_2 \rightarrow A_3, \eta:
A_3 \rightarrow A_4,$
 be proper norm-symmetric
embeddings of 4-cycle algebras.  If the compositions are
$\epsilon_{\ast}$-close to generic embeddings then either 
$\eta$ is $\epsilon$-strict
or $\phi$ has even multiplicity.  
\end{lma}

\begin{pf}
Suppose that the compositions are $\epsilon_*$-close to generic embeddings
and that $\eta$ is not $\epsilon$-strict.
The composition
$\eta \circ \psi$ is $\epsilon_*$-strict and $\epsilon_* \le 
\epsilon^2$ and so by Lemma 6.3 the map $\psi$ is not $\epsilon_1$-strict
where $\epsilon_1 = (\epsilon /51)^2$.
Since $\psi \circ \phi$ is $\epsilon_*$-strict it 
follows similarly, since $\epsilon_* \le   
\epsilon_1^2$ and $\psi$ is not 
$\epsilon_1$-strict, that $\phi$ is not $\epsilon_2$-strict
where $\epsilon_2 = (\epsilon_1 /51)^2$.

We can now continue in a similar fashion to the proof of Lemma 3.1.
Adopting the notation there, since $\psi$ is $\epsilon_1$-strict
we may assume that $a_1$ is not $\epsilon_1$-close to a partial isometry.
On the other hand our hypotheses imply that the matrix
\[
v_1 = \matrrcc{\alpha_1 \otimes a_1}{\beta_1 \otimes b_1}{\delta_1
\otimes d_1}{\gamma_1 \otimes c_1}
\]
is $\epsilon_*$-close to a partial isometry. We now show, as before,
that $\alpha_1$ is a strict contraction.

Suppose that this is not the case.
Then, as before, the matrix $a_1$ appears as an orthogonal part
of $v_1$. That is, there are projections $p, q$ such that
\[
v_1 = qv_1p + (1-q)v_1(1-p)
\]
and $qv_1p = a_1$.
However, if
\[
 \matrrcc
{z_1}{z_2}{z_3}{z_4}
\]
is a partial isometry which is $\epsilon_*$-close to 
\[
v_1 = \matrrcc     
{qv_1p}{0}{0}{(1-q)v_1(1-p)}
\]
then $\|z_2\| \le \epsilon_*$ and $\|z_3\| \le \epsilon_*$,
and so Lemma 6.8 implies
that $z_1$ is $16\epsilon_*$-close to a partial isometry.
Thus $ qv_1p$ is $17\epsilon_*$-close to a partial isometry, contrary 
to the fact that $17\epsilon_* < \epsilon_1$.

Similarly, considering $w_1, x_1, y_1$, it follows that
$\alpha_2, \alpha_3, \alpha_4$
are strict contractions. By the properness of 
$\phi$ each $\alpha_i$
is nonzero and so each of the matrices of Table 1 is nonzero and none is a partial
isometry. Thus, by the argument of Lemma 3.1,
$\phi$ has even multiplicity.
\end{pf}

We also require an approximate version of Lemma 3.2. For the proof of
this the 
following three approximation principles will be convenient.

First, recall that if $e,f$ are projections in a $C^{\ast}$-algebra with
$||e-f||<1$ then there is a partial isometry $v$ in the algebra with  
initial projection $e$, final projection $f$ and $||v-e||<2||e-f||$. 
From this it follows that if $\pi_1$ and $\pi_2$
are star-homomorphisms between finite-dimensional
C*-algebras which are $\epsilon$-close  then there is a
unitary $u$, with 
$\|1-u\| < 2 $, such that 
$\pi_1(a) = u\pi_2(a)u^*$ for all $a$ .

Second, if a contraction $v$ is $\epsilon$-close to a partial isometry
$z$ then the range projection of $v$, $rp(v)$, is $(2\epsilon)^{1/2}$-close to
$zz^*$.

Finally, we need the following elementary lemma.
\medskip

\begin{lma}
Let $E_1, \dots ,E_4$ be projections for which
$\|E_1 + E_2 - E_3 - E_4\| \le 2\epsilon$
and $\|E_1E_4\| \le \epsilon,
\|E_2E_3\| \le \epsilon.$
Then $\|E_1 - E_3\| \le (6\epsilon)^{1/2}$ \ \ and \ \ 
$\|E_2 - E_4\| \le (6\epsilon)^{1/2}$.
\end{lma}
\medskip

\begin{pf}
Let $x$ be a unit vector with $E_2x = x$, so that $\|E_3x\| \le \epsilon.$
Then
\[
2\epsilon \ge \langle(E_1 + E_2 - E_3 - E_4)x,x\rangle
\]
and so
\[
3\epsilon \ge \langle(E_1 + E_2 - E_4)x,x\rangle 
\ge \langle(E_2 - E_4)x,x\rangle .
\]
Thus
\[
\langle(E_2 - E_4E_2)z,z\rangle \le 3\epsilon
\]
for all unit vector $z$.

Similarly $\langle(E_2 - E_2E_4)z,z\rangle \le 3\epsilon$ and so
$\langle(E_2 - E_4)^2z,z\rangle \le 6\epsilon.$
\end{pf}
\medskip

\begin{lma}
Let $\phi: A_1 \rightarrow A_2, \psi: A_2
\rightarrow A_3$ be proper $\epsilon$-strict star extendible embeddings between
4-cycle algebras.  If the composition $\psi \circ \phi$ 
is $ \epsilon$-close to
a generic star extendible embedding then $\psi$ is $g(\epsilon)$-close
to a rigid embedding 
where $g(t)$ is a nonnegative continuous function on $[0,1]$ with
$g(0) = 0$.
\end{lma}

\begin{pf}
First  assume that $A_1 = A(D_4)$ and that the
maps are unital.  Using the usual notation of Section 3, 
since the composition $\psi \circ \phi$ is $
\epsilon$-close to a rigid embedding it follows that the pair 

$$v = \left[ \begin{array}{cc} v_1 & v_2 \\ v_4 & v_3 \end{array}
\right], \quad \quad w = \left[ \begin{array}{cc} w_1 & w_2 \\
w_4 & w_3 \end{array} \right] $$
is $\epsilon$-close to a pair

$$\tilde{v} = \left[ \begin{array}{cc} \tilde{v}_1 & \tilde{v}_2 \\ 
\tilde{v}_4 & \tilde{v}_3 \end{array}
\right], \quad \quad \tilde{w} = \left[ \begin{array}{cc} 
\tilde{w}_1 & \tilde{w}_2 \\
\tilde{w}_4 & \tilde{w}_3 \end{array} \right]$$
where $\tilde{v} = \lambda (e_{13}), \tilde{w} = 
\lambda (e_{14})$ and $\lambda$ is a
rigid embedding.  Moreover, by the remarks preceding the lemma, 
at the
expense of replacing $ \epsilon$ by $4 \epsilon$ 
we may assume that 
$\lambda$ and $\psi \circ \phi$ agree on $A \cap A^*$.
In particular
 $v, \tilde{v}$ and $w, \tilde{w}$ have the same initial 
projections and the same
final projections.  

Since $\lambda$ is rigid the final projection of
$\tilde{v}_1$ is equal to the final projection of $\tilde{w}_2$.
On the other hand since $v_1$ is $4\epsilon$-close to the partial isometry
$\tilde{v_1}$
the range projection \ $rp(v_1)$ is $K_1\epsilon^{1/2}$-close
to the final projection of $\tilde{v}_1$. ($K_1 = 3$ will do.)
Thus the range projections
\[
rp(v_1)  
= rp \matrrcc{\alpha_1 \otimes a_1}{\beta_1 \otimes b_1}{\delta_1  
\otimes d_1}{\gamma_1 \otimes c_1},\ \ \  
rp(w_2) = rp \matrrcc{\alpha_2 \otimes a_2}{\beta_2 \otimes b_2}{\delta_2
\otimes d_2}{\gamma_2 \otimes c_2},\ \ \
\]
are $2K_1 \epsilon^{1/2}$-close.

On the other hand, the hypotheses imply that
 all of the matrices of Table 1 and Table 2 are
$\epsilon$-close to partial isometries, say $||\alpha_i - \alpha'_i || \leq 
\epsilon,  \dots, || \delta_i - \delta'_i || \leq \epsilon , 
||a_i - a_i'|| \leq \epsilon , \dots,
||b_i - b_i'|| \leq \epsilon .$
Thus, with the obvious notation,
$\|v_1' - v_1\| \le 4\epsilon$
and $\|w_1' - w_1\| \le 4\epsilon$
and so $\|v_1' - \tilde{v}_1\| \le 5\epsilon, 
\|w_1' - \tilde{w}_1\| \le 5\epsilon$ and 
\[
\|rp(v_1') - rp(\tilde{v}_1)\| \le (10\epsilon)^{1/2},\ \ 
\|rp(w_1') - rp(\tilde{w}_1)\| \le (10\epsilon)^{1/2}.
\]
Thus
\[
rp \matrrcc{\alpha_1' \otimes a_1'}{\beta_1' \otimes b_1'}{\delta_1'
\otimes d_1'}{\gamma_1' \otimes c_1'},\ \ \
 rp \matrrcc{\alpha_2' \otimes a_2'}{\beta_2' \otimes b_2'}{\delta_2'
\otimes d_2'}{\gamma_2' \otimes c_2'}\ \ \
\]
are $2(10\epsilon)^{1/2}$-close.
But $V_1'$, which need not be a contraction, is 
$2\epsilon$-close to a partial isometry and so its range projection
must be close to 
$v_1'{v_1'}^*$. In fact calculation shows that 
it is $7\epsilon^{1/2}$-close. 
It follows then that
\[
\|v_1'{v_1'}^* - w_2'{w_2'}^*\| \le 7\epsilon^{1/2} + 2(10\epsilon)^{1/2}
+ 7\epsilon^{1/2} \le 22\epsilon^{1/2}.
\]
In particular
\[
\|(\alpha_1'{\alpha_1'} ^* \otimes a_1'{a_1'} ^* +
\beta_1'{\beta_1'} ^* \otimes b_1{'b_1'} ^*) - 
(\alpha_2'{\alpha_2'} ^* \otimes a_2'{a_2'} ^* +
\beta_2'{\beta_2'} ^* \otimes b_2'{b_2'} ^*)\| \le 22\epsilon^{1/2}.
\]
However $a_1', \dots a_4'$ are partial isometries which are the entries
of a matrix which is $2\epsilon$-close to a partial isometry
and so $a_1$ and $a_2$ have almost orthogonal final projections.
In fact
$\|q_1q_2\| \le (10\epsilon)^{1/2}$ where $q_i = rp(a_i'), i=1,2.$ Likewise
$\|e_1e_2\| \le (10\epsilon)^{1/2}$ where $e_i = rp(b_i'), i=1,2.$
From the lemma above
it follows that
\[
\|(\alpha_1'{\alpha_1'} ^* \otimes a_1'{a_1' }^* +
\beta_2'{\beta_2' }^* \otimes b_2'{b_2' }^*)\| \le (3(22\epsilon)^{1/2})^{1/2}
\le K\epsilon^{1/4}.
\]

In a similar way it follows that  $b_2', c_3'$ 
have $K\epsilon^{1/4}$-close initial projections,
and $c_3', d_4'$ have close  final projections.
We may thus rechoose $b_2', c_3',$ so that these initial and final projections
match. Set $d_4' = c_3'(b_2')^*a_1'$ to obtain
a 4-cycle of partial isometries $a_1', b_2', c_3', d_4'$
approximating $a_1, b_2, c_3, d_4$. This 4-cycle of
partial isometries defines a rigid embedding $\zeta_1$ 
(of signature $\{\mbox{rank}(a_1'),0,0,0\}$) which is 
an approximate summand of
$\psi$. 

Similarly we may construct an approximating 
4-cycle $b_1',a_2',d_3', c_4'$ and an associated embedding
$\zeta_2$ and we can arrange, moreover, that $\zeta_2, \zeta_1$
are orthogonal, so that $\zeta_1 \oplus \zeta_2$
is a rigid embedding.
Continuing
in this way, obtain the approximating rigid embedding
$\zeta_1 \oplus \zeta_2 \oplus \zeta_3 \oplus \zeta_4$.
At each stage the approximation factors deteriorate by at most a 
multiplicative factor or by a fourth root and so the lemma
is obtained with $g(t) = Ct^{1/n}$ for some constant $C$
and some integer $n$.
\end{pf}
\medskip

 
\begin{thm}
Let $\ca, \ca{'}$ be odd generic systems of 4-cycle
algebras with direct limits
$A = \indlimit \ca \ , \  A' =
\indlimit \ca{'}$
respectively.  Then
$A$ and $A{'}$ are star extendibly
isomorphic if and only if
$\ca$ and $\ca{'}$ are regularly
isomorphic
\end{thm}
\medskip

Whilst it is possible to use the lemmas 
above to obtain a direct proof of the theorem we prefer
to show the key fact that
star extendibly isomorphic limit algebras
have isomorphic $K_0H_1$ groups and scales.
In particular these groups and
scales are indeed invariants for the algebra.
In fact  the arguments given show that close limit algebras of this type 
have isomorphic $K_0H_1$ invariants and so this approach leads 
to the fact that close odd generic limit algebras are isomorphic.

Let $x$ be an operator on ${\Bbb C}^n$ with $\|x - u\| \le \frac{1}{3}$
for some partial isometry $u$.
Then the rank of $u$ can be recovered from $x$ as the rank
of the spectral projection for $|x|$ for 
the interval $(\frac{1}{2}, \frac{3}{2})$.
Refer to this as the {\it upper rank}
of $x$, denoted $rank_u(x)$.
From this simple observation it follows that if $\phi : A_1 \to A_2$ is
a star extendible embedding between 4-cycle algebras which
is $\frac{1}{3}$-close to a rigid embedding  $\phi '$ then,
in the usual notation,
\[
H_1(\phi ') = [rank_u(\alpha_1) - rank_u(\beta_1) + rank_u(\gamma_1)
- rank_u(\delta_1)].
\]
Thus generic systems $\ca, \ca'$ which are linked by an approximately commuting
diagram, with approximately rigid crossover
maps, have isomorphic $H_1$ groups. 
Since close maps induce the same
$K_0$ maps it follows that the jointly scaled $K_0H_1$ invariants for 
$\ca, \ca'$ are isomorphic and so, by
Theorem 5.2, $\ca, \ca'$ are star extendibly isomorphic.
Note that the  approximately commuting
diagram need not be asymptotically commuting; $\frac{1}{3}$-close 
approximations are 
all that are required.
\medskip

{\it Proof of the Theorem 6.7}\ \ 
It will be enough to obtain an approximately commuting
diagram for $\ca, \ca'$.
More precisely, it will suffice to obtain
embeddings $\eta_1, \eta_2, \dots $ with
\[
A_{t_1} @>{\eta_1}>> A_{s_1}' @>{\eta_2}>>  
        A_{t_2} @>{\eta_3}>> A_{s_2}' \to \dots
\]
such that each $\eta_k$ is 
$\frac{1}{3}$-close to a rigid embedding and each 
composition $\eta_{k+1} \circ \eta_k$ is $\frac{1}{3}$-close
to the given generic embedding.

Let $\theta : A \to A'$  be a star extendible isomorphism.
Note that $A$ admits a direct sum decomposition $A = A \cap A^* + A^r$ 
where $A^r$ is the largest ideal whose square is trivial. The ideal $A^r$ is also the 
Jacobson radical. Thus $\theta$ maps $A \cap A^*$ to $A' \cap {A'}^*$  and
$\theta (A^r) = {A'}^r$. Since the restriction  $\theta : A_1 \to A'$
maps $A_1^r$ into $(A')^r$ it is elementary to construct in
$(A_k')^r$, for some suitably large $k$, partial isometries $X, Y, Z$ 
which are $\epsilon_1$-close to
$\theta (e_{13}), \theta (e_{14}), \theta (e_{24})$.
Note that $W = ZY^*X$ also belongs to $(A_k')^r$, and so there is a star extendible proper embedding
$\phi_1 : A_1 \to A_{k}{'}^r$ which is $K\epsilon$-close to the restriction
of $\theta$ to $A_1$.
Composing $\phi_1$ with the given generic embedding $A_k' \to A_{k+1}'$
we may also assume that $\phi$ is norm-symmetric.
Repeating this idea twice more we obtain
the double triangle
\[
A_{1} @>{\phi}>> A_{n_1}' @>{\psi}>>
        A_{m_1} @>{\eta}>> A_{n_2}' 
\]
with $\phi, \psi, \eta$ proper, norm-symmetric,
star extendible
and with
the compositions $\psi \circ \phi, \eta \circ \psi$ close to the given
generic embeddings.

Let $\epsilon > 0$, with  $\epsilon < \frac{1}{16}$. For suitably small $\epsilon_1$
it follows from Lemma 6.4 that $\eta$ is $\epsilon$-strict.

In this way construct $\epsilon$-strict proper star extendible
embeddings $\eta_1 : A_{t_1} \to A_{s_1}',\  \eta_2 : A_{s_1}' \to
A_{t_2}, \dots,$
for which the compositions 
$\eta_{k+1} \circ \eta_k$
are $\epsilon$-close to the given rigid embeddings. By Lemma
6.6  the embeddings $\eta_k$ are $g(\epsilon)$-close to
rigid embeddings. For  $\epsilon$
 sufficiently small the maps $\eta_k$ give the desired approximately commuting diagram.
\hfill $\Box$.
\medskip

Combining the above with Theorem 5.2 one  obtains a 
classification of the operator algebras of odd systems.
In particular, in the unital case we have
\medskip

\begin{thm}
Odd generic unital systems of 4-cycle algebras have operator
algebra limits which are star extendibly isomorphic
if and only if their invariants
\[
(K_0(-) \oplus H_1(-), \Sigma_0(-), \Sigma(-))
\]
are  isomorphic.
\end{thm}
\medskip

A more detailed perturbational analysis 
should lead to a full generalisation of 
Lemma 3.1 from which the even case would similarly follow.
This in turn will enable the extension of the relative position 
analysis of sections 4 and 5 to matroid C*-algebras.

\end{document}